\documentclass[12pt]{article}
\usepackage[top=1in,bottom=1in,left=1in,right=1in]{geometry}
\usepackage{indentfirst}
\usepackage{amsfonts,amsmath,amsthm,amssymb,hyperref,enumerate,bm}
\usepackage{longtable}

\usepackage{color}

\newtheorem{lemma}{Lemma}[section]
\newtheorem{theorem}{Theorem}[section]

\newtheorem{proposition}{Proposition}[section]

\theoremstyle{definition}
\newtheorem{construction}{Construction}[section]
\newtheorem{hypothesis}{Hypothesis}[section]

\let\leq=\leqslant 
\let\geq=\geqslant

\newcommand{\Soc}{{\rm Soc}}

\newcommand{\B}{\mathcal{B}}
\newcommand{\D}{\mathcal{D}}
\newcommand{\C}{\mathcal{C}}
\def\P{\mathcal{P}}

\newcommand{\PG}{{\rm PG}}
\def\proof{\par\noindent{\sc Proof.~}}
\newcommand{\Aut}{{\rm Aut}}

\newcommand{\lcm}{{\rm lcm}}

\newcommand{\PSL}{{\rm PSL}}
\newcommand{\PGL}{{\rm PGL}}
\newcommand{\PGaL}{{\rm P\Gamma L}}

\newcommand{\SL}{{\rm SL}}

\newcommand{\GL}{{\rm GL}}

\newcommand{\GGL}{{\rm \Gamma L}}

\newcommand{\AGGL}{{\rm A\Gamma L}}

\newcommand{\Sp}{{\rm Sp}}
\newcommand{\normal}{\trianglelefteq}
\def\gbc{\genfrac{[}{]}{0pt}{}}
\def\la{\langle}
\def\ra{\rangle}

\DeclareMathOperator{\Out}{Out}
\begin{document}

\title{On flag-transitive $2$-$(v,k,2)$ designs}

\author{Alice Devillers\thanks{The first and third author were supported by an ARC Discovery Grant Project DP200100080.}
            \and Hongxue  Liang\thanks{The second author was supported by the NSFC (Grant No.11871150).} 
           \and Cheryl E. Praeger \and Binzhou Xia
\date{}
}

\newcommand{\Addresses}{{
  \bigskip
  \footnotesize

  Devillers, Praeger: \textsc{Centre for Mathematics of Symmetry and Computation,
University of Western Australia,
35 Stirling Highway,
Perth 6009, Australia.}\par\nopagebreak
   Email: {\tt \{Alice.Devillers,\;Cheryl.Praeger\}@uwa.edu.au}

  \medskip

  Liang (Corresponding author): \textsc{School of Mathematics and Big Data, Foshan University, Foshan 528000, P. R. China}\par\nopagebreak
  Email: {\tt hongxueliang@fosu.edu.cn}

  \medskip

Xia: \textsc{School of Mathematics and Statistics, The University of Melbourne, Parkville, VIC 3010, Australia}\par\nopagebreak
   Email: {\tt binzhoux@unimelb.edu.au}

}}

\maketitle

\begin{abstract}
This paper is devoted to the classification of flag-transitive $2$-$(v,k,2)$ designs.
We show that apart from two known symmetric $2$-$(16,6,2)$ designs, every flag-transitive subgroup $G$ of the automorphism group of a nontrivial $2$-$(v,k,2)$ design is primitive of affine or almost simple type. Moreover, we classify the $2$-$(v,k,2)$ designs admitting a flag transitive almost simple group $G$ with socle $\PSL(n,q)$ for some $n\geq 3$. Alongside this analysis we give a construction for   a flag-transitive $2$-$(v, k-1, k-2)$ design from a given flag-transitive $2$-$(v,k, 1)$ design which induces a 2-transitive action on a line. Taking the design of points and lines of the projective space $\PG(n-1,3)$ as input to this construction yields a
$G$-flag-transitive $2$-$(v,3,2)$ design where $G$ has socle $\PSL(n,3)$ and $v=(3^n-1)/2$. Apart from these designs, our $\PSL$-classification yields exactly one other example, namely the complement of the Fano plane.



\end{abstract}

\section{Introduction}

A \emph{$2$-$(v,k,\lambda)$ design} $\D$ is a pair $(\P,\B)$ with a set $\P$ of $v$ \emph{points} and a set $\B$ of \emph{blocks}
such that each block is a $k$-subset of $\P$ and each two distinct points are contained in $\lambda$ blocks.
We say $\D$ is \emph{nontrivial} if $2<k<v$, and \emph{symmetric} if $v=b$.
All $2$-$(v,k,\lambda)$ designs in this paper are assumed to be nontrivial.
An automorphism of $\D$ is a permutation of the point set
which preserves the block set.
The set of all automorphisms of $\D$ with the composition of permutations forms a
group, denoted by $\Aut(\D)$.
For a subgroup $G$ of $\Aut(\D)$,
$G$ is said to be \emph{point-primitive}
if $G$ acts primitively on $\P$,
and said to be \emph{point-imprimitive} otherwise.
A \emph{flag} of $\D$ is a point-block
pair $(\alpha,B)$ where $\alpha$ is a point
and $B$ is a block incident with $\alpha$.
A subgroup $G$ of $\Aut(\D)$ is said to be \emph{flag-transitive} if $G$ acts transitively on the set of flags of $\D$.

A $2$-$(v,k,\lambda)$ design with $\lambda=1$
is also called a finite \emph{linear space}.
In 1990, Buekenhout, Delandtsheer, Doyen, Kleidman, Liebeck and Saxl~\cite{BDDKLS} classified all flag-transitive linear spaces apart from
those with a one-dimensional affine automorphism group.
Since then, there have been efforts to classify $2$-$(v,k,2)$ designs $\mathcal{D}$
admitting a flag-transitive group $G$ of automorphisms.
Through a series of papers~\cite{Re2005,Biplane1,Biplane2,Biplane3}, Regueiro proved that,
if $\mathcal{D}$ is symmetric, then either $(v,k)\in\{(7,4),(11,5),(16,6)\}$, or $G\leq \AGGL(1,q)$ for some odd prime power $q$. Recently, Zhou and the second author~\cite{Liang1} proved that, if $\mathcal{D}$
is not symmetric and $G$ is point-primitive, then $G$ is affine or almost simple. In each of these cases $G$ has a unique minimal normal subgroup, its \emph{socle} $\Soc(G)$, which is elementary abelian or a nonabelian simple group, respectively.

Our first objective in this paper is to fill in a missing piece in this story, namely to treat the case where $G$ is flag-transitive and point-imprimitive and $\mathcal{D}$ is a not-necessarily-symmetric $2$-$(v,k,2)$ design. Such flag-transitive, point-imprimitive designs exist: it was shown in 1945 by Hussain \cite{Huss}, and independently in 1946 by Nandi \cite{Nandi},
that there are exactly three $2$-$(16,6,2)$-designs. O'Reilly Regueiro
\cite[Examples 1.2]{Re2005} showed that exactly two of these designs are flag-transitive, and each admits a point-imprimitive, flag-transitive subgroup of automorphisms (one with automorphism group $2^4\mathrm{S}_6$ and point stabiliser $(\mathbb{Z}_2\times\mathbb{Z}_8)(\mathrm{S}_4.2)$ and the other with automorphism group $\mathrm{S}_6$ and point stabiliser $\mathrm{S}_4.2$, see also \cite[Remark 1.4(1)]{Praeger}).
We prove that these are the only point-imprimitive examples, and thus, together with~\cite[Theorem 1.1]{Liang1} and~\cite[Theorem 2]{Re2005}, we obtain the following result.

\begin{theorem}\label{Th1}
Let $\mathcal{D}$ be a $2$-$(v,k,2)$ design with a flag-transitive group $G$ of automorphisms. Then  either
\begin{enumerate}
\item[\rm(i)] $\mathcal{D}$ is one of two known symmetric $2$-$(16,6,2)$ designs with $G$ point-imprimitive; or
\item[\rm(ii)] $G$ is point-primitive of affine or almost simple type.
\end{enumerate}
\end{theorem}

Theorem~\ref{Th1} reduces the study of flag-transitive $2$-$(v,k,2)$ designs to those whose automorphism group $G$ is point-primitive of affine or almost simple type. Regueiro \cite{Re2005,Biplane1,Biplane2,Biplane3} has classified all such examples where the design is symmetric (up to those admitting a one-dimensional affine group). In the non-symmetric case, the second author and Zhou have dealt with the cases where the socle $\Soc(G)$ is a sporadic simple group or an alternating group, identifying three possibilities: namely $(v,k)=(176,8)$ with $G=\mathrm{HS}$, the Higman-Sims group in \cite{Liang1}, and $(v,k)=(6,3)$ or $(10,4)$ with $\Soc(G)=A_v$ in~\cite{Liang2}. Our contribution is the case where $\Soc(G)=\PSL(n,q)$ for some $n\geq 3$ and $q$ a prime power.
In contrast to the cases considered previously,
an infinite family of examples occurs, which may be obtained from the following
general construction method for flag-transitive designs from linear spaces.

\begin{construction}\label{cons}
For a $2$-$(v,k,1)$ design $\mathcal{S=(P,L)}$ with $k\geq 3$, let
\[
\mathcal{B}=\{\ell\setminus\{\alpha\}\,\mid\,\ell\in \mathcal{L},\,\alpha\in\ell\}
\]
and $\mathcal{D(S)}=(\P,\B)$.
\end{construction}

We show in Proposition \ref{exist 1}
that $\mathcal{D(S)}$ is a $2$-$(v,k-1,k-2)$ design, and
moreover, that $\mathcal{D(S)}$ is $G$-flag-transitive whenever $G\leq\Aut(\mathcal{S})$ is flag-transitive on $\mathcal{S}$ and induces a 2-transitive action on each line of $\mathcal{S}$. In particular, these conditions hold if $\cal S$ is the design of points and lines of $\PG(n-1,3)$,  for some $n\geq 3$, and $\Soc(G)=\PSL(n,3)$   (Proposition~\ref{exist 1}). Apart from these designs, our analysis shows that there is only one other $G$-flag-transitive $2-(v,k,2)$ design with $\Soc(G)=\PSL(n,q)$, $n\geq3$.

\begin{theorem}\label{Th2}
Let $\D$ be a $2$-$(v,k,2)$ design admitting a flag-transitive group $G$ of automorphisms, such that $\Soc(G)=\PSL(n,q)$ for some $n\geq3$ and prime power $q$. Then either
\begin{enumerate}[(a)]
\item $\D=\D(\cal S)$ is as in Construction~$\ref{cons}$, where $\cal S$ is the design of points and lines of $\PG(n-1,3)$; or
\item $\D$ is the complement of the Fano plane (that is, blocks are the complements of the lines of $\PG(2,2)$).
\end{enumerate}
\end{theorem}
The designs in part (a) are non-symmetric (Proposition~\ref{exist 1}), while the complement of the Fano plane is symmetric, and arises also in Regueiro's classification \cite[Theorem 1]{Biplane2} (noting that the group $\PSL(3,2)$ is isomorphic to the group $\PSL(2,7)$ in her result).




The proofs of Theorems \ref{Th1} and \ref{Th2} will be given in Sections \ref{sec3} and \ref{sec4}, respectively.

\section{Preliminaries}\label{sec2}

We first collect some useful results on flag-transitive designs and groups of Lie type.

\begin{lemma} \label{condition 1}
Let $\D$ be a $2$-$(v,k,\lambda)$ design and let $b$ be the number of blocks of $\D$.
Then the number of blocks containing each point of $\D$ is a constant $r$ satisfying the following:
\begin{enumerate}
\item[\rm(i)] $r(k-1)=\lambda(v-1)$;
\item[\rm(ii)] $bk=vr$;
\item[\rm(iii)] $b\geq v$ and $r\geq k$;
\item[\rm(iv)] $r^2>\lambda v$.
\end{enumerate}
In particular, if $\D$ is non-symmetric then $b>v$ and $r>k$.
\end{lemma}

\proof
Parts~(i) and~(ii) follow immediately by simple counting.
Part~(iii) is  Fisher's Inequality \cite[p.99]{Ryser}.
By~(i) and~(iii) we have
\[
r(r-1)\geq r(k-1)=\lambda(v-1)
\]
and so $r^2\geq\lambda v+r-\lambda$.
Since $\D$ is nontrivial, we deduce from (i) that $r>\lambda$.
Hence $r^2>\lambda v$, as stated in part~(iv).
\qed

For a permutation group $G$ on a set $\P$ and an element $\alpha$ of $\P$, denote by $G_\alpha$ the stabiliser of $\alpha$ in $G$, that is, the subgroup of $G$ fixing $\alpha$. A \emph{subdegree} $s$ of a transitive permutation group $G$ is the length of some orbit of $G_\alpha$. We say that $s$ is \emph{non-trivial} if the orbit is not $\{\alpha\}$, and $s$ is \emph{unique} if $G_\alpha$ has only one orbit of size $s$.

\begin{lemma} \label{condition 2}
Let $\D$ be a $2$-$(v,k,\lambda)$ design, let $G$ be a flag-transitive subgroup of $\Aut(\D)$, and let $\alpha$ be a point of $\D$.
Then the following statements hold:
\begin{enumerate}
\item[\rm(i)] $|G_\alpha|^3>\lambda |G|$;
\item[\rm(ii)] $r$ divides $\gcd(\lambda(v-1),|G_{\alpha}|)$;
\item[\rm(iii)] $r$ divides $\lambda\gcd(v-1,|G_{\alpha}|)$;
\item[\rm(iv)] $r$ divides $s\gcd(r,\lambda)$ for every nontrivial subdegree $s$ of $G$.
\end{enumerate}
\end{lemma}

\proof
By Lemma~\ref{condition 1} we have $r^2>\lambda v$.
Moreover, the flag-transitivity of $G$ implies that $v=|G|/|G_\alpha|$ and $r$ divides $|G_\alpha|$, and in particular, $|G_\alpha|\geq r$.
It follows that
\[
|G_\alpha|^2\geq r^2>\lambda v=\frac{\lambda|G|}{|G_\alpha|}
\]
and so $|G_\alpha|^3>\lambda|G|$.
This proves statement~(i).

Since $r$ divides $r(k-1)=\lambda(v-1)$ and $r$ divides $|G_\alpha|$, we conclude that
$r$ divides
\begin{equation}\label{1}
\gcd(\lambda(v-1),|G_{\alpha}|),
\end{equation}
as statement~(ii) asserts. Note that the quantity in~\eqref{1} divides
\[
\gcd(\lambda(v-1),\lambda|G_{\alpha}|)=\lambda\gcd(v-1,|G_{\alpha}|).
\]
We then conclude that $r$ divides $\lambda\gcd(v-1,|G_\alpha|)$, proving statement~(iii).

Finally, statement~(iv) is proved in \cite[p.91]{Dav1} and \cite{Dav2}.
\qed

For a positive integer $n$ and prime number $p$, let $n_p$ denote the \emph{$p$-part of $n$} and let $n_{p'}$ denote the \emph{$p'$-part of $n$}, that is, $n_p=p^t$ such that $p^t\mid n$ but $p^{t+1}\nmid n$ and $n_{p'}=n/n_p$.
We will denote by $d$  the greatest common divisor of $n$ and $q-1$.

\begin{lemma}\label{bound}
Suppose that $\D$ is a $2$-$(v,k,2)$ design admitting a flag-transitive point-primitive group $G$ of automorphisms with socle $X=\PSL(n,q)$, where $n\geq 3$ and $q=p^f$ for some prime $p$ and positive integer $f$, and $d=\gcd(n,q-1)$. Then for any point $\alpha$ of $\D$ the following statements hold:
\begin{enumerate}
\item[\rm(i)] $|X|<2(df)^2|X_{\alpha}|^3$;
\item[\rm(ii)] $r$ divides $2df|X_{\alpha}|$;
\item[\rm(iii)] if $p\mid v$, then  $r_p$ divides $2$, $r$ divides $2df|X_{\alpha}|_{p'}$, and $|X|<2(df)^2|X_{\alpha}|^2_{p'}|X_{\alpha}|$.
\end{enumerate}
\end{lemma}

\proof
Since $G$ is point-primitive and $X$ is normal in $G$, the group $X$ is transitive on the point set. Hence $G=XG_\alpha$ and so
\[
\frac{|G_\alpha|}{|X_\alpha|}=\frac{|G_\alpha|}{|X\cap G_\alpha|}=\frac{|XG_\alpha|}{|X|}=\frac{|G|}{|X|}.
\]
Moreover, as $\Soc(G)=X=\PSL(n,q)$, we have $G\leq\Aut(X)$. Hence $|G_\alpha|/|X_\alpha|=|G|/|X|$ divides $|\Out(X)|=2df$.
Consequently, $|G_\alpha|/|X_\alpha|\leq2df$.
Since Lemma \ref{condition 2}(i) yields
\[
|G_\alpha|^3>2|G|=\frac{2|X||G_\alpha|}{|X_\alpha|},
\]
it follows that
\[
2|X|<|X_\alpha||G_\alpha|^2=\left(\frac{|G_\alpha|}{|X_\alpha|}\right)^2|X_\alpha|^3\leq (2df)^2|X_\alpha|^3.
\]
This leads to statement~(i).
Since $|G_\alpha|/|X_\alpha|$ divides $|\Out(X)|=2df$ and the flag-transitivity of $G$ implies that $r$ divides $|G_\alpha|$, we derive that $r$ divides $2df|X_{\alpha}|$, as in statement~(ii).

Now suppose that $p$ divides $v$.
Then the equality $2(v-1)=r(k-1)$ implies that $r_p$ divides $2$.
As a consequence of this and part (ii) we see that
$r$ divides $2df|X_{\alpha}|_{p'}$.
Since $r^2>2v$ by Lemma~\ref{condition 1}(iv), and $v=|X|/|X_{\alpha}|$ by the point-transitivity of $X$, it then follows that
\[
(2df|X_{\alpha}|_{p'})^2>2v=\frac{2|X|}{|X_\alpha|}.
\]
This implies that $2(df)^2|X_{\alpha}|^2_{p'}|X_{\alpha}|>|X|$, completing the proof of part (iii).
\qed

\begin{lemma}\label{L:subgroupdiv}
Suppose that $\D$ is a $2$-$(v,k,2)$ design admitting a flag-transitive point-primitive group $G$ of automorphisms with socle $X=\PSL(n,q)$, where $n\geq 3$ and $q=p^f$ for some prime $p$ and positive integer $f$, and $d=\gcd(n,q-1)$. Let $\alpha$ and $\beta$ be distinct points of $\D$, and suppose $H\leq G_{\alpha,\beta}$. Then $r$ divides $4df|X_\alpha|/|H|$.
\end{lemma}
\proof By Lemma~\ref{condition 2}(iv), $r$ divides $2|\beta^{G_\alpha}|=2|G_\alpha|/|G_{\alpha\beta}|$. Since
$|{G}_\alpha|$ divides $2df|{X}_\alpha|$ (see proof of Lemma \ref{bound}) and $H$ divides $|G_{\alpha,\beta}|$,
it follows that $r$ divides $4df|X_\alpha|/|H|$.
\qed
%

We will need the following results on finite groups of Lie type.

\begin{lemma}\label{parabolic} Suppose that $\D$ is a $2$-$(v,k,2)$ design admitting a flag-transitive point-primitive group $G$ of automorphisms with socle $X=\PSL(n,q)$, where $n\geq 3$ and $q=p^f$ for some prime $p$ and positive integer $f$, and $r$ is the number of blocks incident with a given point. Let $\alpha$ be a point of $\D$.
Suppose that $X_\alpha$ has a normal subgroup $Y$, which is a finite simple group of
Lie type in characteristic $p$, and $Y$ is not isomorphic to $\mathrm{A}_{5}$ or $\mathrm{A}_{6}$ if $p=2$.
If $r_p\mid 2_p$, then $r$ is divisible by the index of
a proper parabolic subgroup of $Y$.
\end{lemma}
\proof
Since $G$ is flag-transitive, we have $r=|G_\alpha|/|G_{\alpha,B}|$, where
$B$ is a block through $\alpha$. Since $X_\alpha\unlhd G_\alpha$, $|X_\alpha|/|X_{\alpha,B}|$ divides $r$.
Now since $Y\unlhd X_\alpha$, we also have that $|Y|/|Y_{B}|$ divides $r$. Let $H:=Y_{B}$. Since $r_{p}\mid 2_{p}$, we have that $|Y{:}H|_{p}\leq 2_{p}$. We claim that $H$ is contained in a proper parabolic subgroup of $Y$.
First assume $|Y{:}H|_{p}=1$. Then by \cite[Lemma 2.3]{Saxl},
$H$ is contained in a proper parabolic subgroup of $Y$. Now suppose $|Y{:}H|_{p}=2$. Then $p=2$ and $4\nmid |Y{:}H|$,
and so by \cite[Lemma 7]{Biplane2}, $H$ is contained in a
proper parabolic subgroup of $Y$. So the claim is proved in both cases.
It follows that $r$ is divisible by the index of a parabolic subgroup of $Y$.
\qed

\begin{lemma}{\rm (\cite[Lemma 4.2, Corollary 4.3]{Alavi})}\label{eq2}
Table~$\ref{tab1}$ gives upper bounds and lower bounds for the orders of certain $n$-dimensional classical groups defined over a field of order $q$, where $n$ satisfies the conditions in the last column.
\begin{table}[h]
\begin{center}
\renewcommand\arraystretch{1.9}
\caption {Bounds for the order of some classical groups}
\label{tab1}
\vspace{5mm}
\begin{tabular}{l|l|l|c}
\hline
Group  $G$& Lower bound on $|G|$  & Upper bound on $|G|$& Conditions on $n$  \\
\hline
$\GL(n,q)$   &$>(1-q^{-1}-q^{-2})q^{n^2}$    &$\leq(1-q^{-1})(1-q^{-2})q^{n^2}$  &$n\geq 2$ \\
$\PSL(n,q)$  &$>q^{n^2-2}$                   &$\leq(1-q^{-2})q^{n^2-1}$          &$n\geq 2$ \\
${\rm GU}(n,q)$   &$\geq(1+q^{-1})(1-q^{-2})q^{n^2}$ &$\leq(1+q^{-1})(1-q^{-2})(1+q^{-3})q^{n^2}$ & $n\geq 2$ \\
${\rm PSU}(n,q)$  &$>(1-q^{-1})q^{n^2-2}$    &$\leq(1-q^{-2})(1+q^{-3})q^{n^2-1}$ &$n\geq 3$ \\
${\rm Sp}(n,q)$&$>(1-q^{-2}-q^{-4})q^{\frac{1}{2}n(n+1)}$  &  $\leq(1-q^{-2})(1-q^{-4})q^{\frac{1}{2}n(n+1)}$ & $n\geq 4$ \\
\hline
\end{tabular}
\end{center}
\end{table}
\end{lemma}

We finish this section with an arithmetic result.
\begin{lemma}{\rm (\cite[Lemma 1.13.5]{Low})}\label{gcd}
Let $p$ be a prime, let $n, e$ and $f$ be positive integers
such that $n>1$ and $e\mid f$,
and let $q_{0}=p^e$ and $q=p^f$.
Then
\begin{enumerate}
\item[\rm(i)]
$\displaystyle{
    \frac{q-1}{\lcm(q_0-1,(q-1)/\gcd(n,q-1))}=\gcd\left(n,\frac{q-1}{q_0-1}\right);
    }$
\item[\rm(ii)]
    $\displaystyle{
    \frac{q+1}{\lcm(q_0+1,(q+1)/\gcd(n,q+1))}=\gcd\left(n,\frac{q+1}{q_0+1}\right);
    }$
\item[\rm(iii)] If $f$ is even, then $q^{1/2}=p^{f/2}$ and \quad
    $\displaystyle{
    \frac{q-1}{\lcm(q^{1/2}+1,(q-1)/\gcd(n,q-1))}=\gcd\left(n,q^{1/2}-1\right).
    }$
\end{enumerate}
\end{lemma}

\section{Proof of Theorem \ref{Th1}}\label{sec3}

Let $\cal D =(\P, \cal B)$ be a $2$-$(v,k,2)$ design admitting
a flag-transitive group $G$ of automorphisms.
If $G$ is point-primitive, then by~\cite{Liang1} and~\cite{Re2005},
$G$ is of affine or almost simple type. Thus we may assume that $G$
leaves invariant a non-trivial partition
$\mathcal{C}=\{\Delta_1,\Delta_2,\dots,\Delta_y\}$ of $\P$,
where
\begin{equation}\label{Eq5}
v=xy.
\end{equation}
with $1<y<v$ and $|\Delta_i|=x$ for each $i$.
If $(v,k)=(16,6)$ then by Lemma~\ref{condition 1}, it follows that
$\cal D$ is symmetric and hence, in the light of the discussion before
the statement of Theorem \ref{Th1}, in this case Theorem \ref{Th1}(i) holds.
Hence we may assume further that  $(v,k)\ne (16,6)$. Our objective now is to derive a contradiction to these assumptions.
Our proof uses the facts, which can easily be verified by \textsc{Magma} \cite{magma}, that for each 2-transitive permutation group of degree $2p=10$ or $22$ there is a unique class of subgroups of index $2p$ and each such group is
almost simple with a 2-transitive unique minimal normal subgroup (its socle). In fact the socle is one of $\PSL(2,9)$ or $\mathrm{A}_{10}$ (for degree 10), or $M_{22}$ or $\mathrm{A}_{22}$ (for degree 22).

First we introduce a new parameter $\ell$: let  $\alpha\in\P$ and $\Delta\in \mathcal{C}$
such that $\alpha\in\Delta$; choose $B\in \mathcal{B}$ containing $\alpha$,
and let $\ell=|B\cap\Delta|$. It follows from~\cite[Lemma~2.1]{Praeger} that,
for each $B'\in\mathcal{B}$ and $\Delta'\in\mathcal{C}$ such that $B'\cap\Delta'\neq\emptyset$,
the intersection size $|B'\cap\Delta'|=\ell$,
so that $B'$ meets each of exactly $k/\ell$ parts of $\cal C$ in $\ell$ points and is disjoint from the other parts. Moreover,
\begin{equation}\label{Eq6}
\ell\mid k\quad \mbox{and}\ 1<\ell < k.
\end{equation}
(Note that the proof of~\cite[Lemma~2.1]{Praeger}
uses flag-transitivity of $\cal D$, but is valid for all $2$-designs, not only symmetric ones.)

\medskip\noindent
\textit{Claim 1:}\quad $(v,b, r,k, \ell)=(x^2, \frac{2x^2(x-1)}{x+2}, 2x-2, x+2, 2)$, and $x=2p$ with $p\in\{5, 11\}$.

\smallskip\noindent
\textit{Proof of Claim:}\quad Counting the point-block pairs $(\alpha',B')$ with $\alpha'\in\Delta\setminus\{\alpha\}$ and $B'$ containing $\alpha$ and $\alpha'$, we obtain
\begin{equation}\label{Eq7}
2(x-1)=r(\ell-1).
\end{equation}
It follows from~\eqref{Eq5} and Lemma~\ref{condition 1}(i) that
\[
r(k-1)=2(xy-1)=2y(x-1)+2(y-1),
\]which together with~\eqref{Eq7} yields
\begin{equation}\label{Eq8}
r(k-1)=yr(\ell-1)+2(y-1).
\end{equation}
Let $z=k-1-y(\ell-1)$. Then $z$ is an integer and, by \eqref{Eq8},  $rz=2(y-1)>0$
so $z$ is a positive integer and
\begin{equation}\label{Eq9}
y=\frac{rz+2}{2}.
\end{equation}
This in conjunction with~\eqref{Eq8} leads to
\[
r(k-1)+2=y(r(\ell-1)+2)=\frac{(rz+2)(r(\ell-1)+2)}{2}.
\]
Hence
\begin{equation}\label{Eq10}
2(k-\ell-z)=rz(\ell-1).
\end{equation}
Since $k\leq r$ (Lemma~\ref{condition 1}(iii)), we have
\[kz(\ell-1)\leq rz(\ell-1)=2(k-\ell-z)<2k,
\]
and hence $z=1$ and $\ell=2$.
Then~\eqref{Eq7} becomes $r=2x-2$,
and so~\eqref{Eq9} gives $y=x$ (and hence $v=x^2$) and the definition of $z$ gives $k=x+2$.
It then follows from $r\geq k$ that $x\geq 4$, and from~\eqref{Eq6} that $k$, and hence also $x$, is even. Finally by Lemma~\ref{condition 1}(ii),
\[
b=\frac{vr}{k}=\frac{x^2(2x-2)}{x+2} = 2x^2-6x+12 -\frac{24}{x+2},
\]
and hence $(x+2)\mid24$. Therefore, $x=4$, $6$, $10$ or $22$, but since we are assuming that $(v,k)\ne (16,6)$ the parameter $x\ne4$. If $x=6$, then $(v,b,r,k)=(36,45,10,8)$, but one can see from~\cite[II.1.35]{Handbook} that there is no $2$-$(36,8,2)$ design. Thus $x=10$ or $22$, and Claim 1 is proved.\qed

\medskip\noindent
\textit{Claim 2:}\quad For $\Delta\in\cal C$, the induced group $G_\Delta^\Delta$ is 2-transitive. Moreover the kernel $K:=G_{(\cal C)}\ne 1$, $\cal C$ is the set of $K$-orbits in $\P$, and  $K^\Delta$ and its socle $\Soc(K)^\Delta$ are $2$-transitive with 2-transitive socle $\PSL(2,9)$ or $\mathrm{A}_{10}$ for degree 10, and  $M_{22}$ or $\mathrm{A}_{22}$ for degree 22. 

\smallskip\noindent
\textit{Proof of Claim:}\quad
Since each element of $G$ fixing $\alpha$ stabilises $\Delta$, we have the inclusion  ${G_\alpha}\leq G_{\Delta}$. Let $\beta, \gamma$ be arbitrary points  in $\Delta\setminus\{\alpha\}$, and consider $B_1\in\mathcal{B}$ containing $\alpha$ and $\beta$, and $B_2\in \mathcal{B}$ containing $\alpha$ and $\gamma$.
Since $G$ is flag-transitive, there exists $h\in G_\alpha$ such that $B_1^h=B_2$, and in particular, $\beta^h\in B_2$. As $\ell=2$ (by Claim 1), each block of $\mathcal{D}$ through $\alpha$ contains exactly one point in $\Delta\setminus\{\alpha\}$. Since $\beta^h\in(\Delta\setminus\{\alpha\})^h=\Delta\setminus\{\alpha\}$, it then follows that $\beta^h=\gamma$. This shows that $G_\alpha$ is transitive on $\Delta\setminus\{\alpha\}$, and hence  $G^{\Delta}_{\Delta}$ is 2-transitive and hence primitive.  

By Claim 1, each non-trivial block of imprimitivity for $G$ in $\P$ has size $x=\sqrt{v}=2p$ (with $p=5$ or $11$),
and hence the induced permutation group $G^{\mathcal{C}}$ on $\mathcal{C}$ is primitive. Suppose  that $K=1$, so $G^{\mathcal{C}}\cong G$.
Since $G$ is point-transitive and $v=4p^2$, it follows that  $|G|=|G^{\mathcal{C}}|$ is divisible by $p^2$, and hence $G^{\mathcal{C}}_{\Delta}\cong G_\Delta$ has order divisible by $p$ (since $|G:G_\Delta|=2p$). Thus $G^{\mathcal{C}}_{\Delta}$ contains an element of order $p$ which acts on $\cal C$ as a $p$-cycle fixing $p$ of the parts. Then by a result of Jordan~\cite[Theorem 13.9]{Wie1964} we have $G^{\mathcal{C}}=\mathrm{A}_{2p}$ or $\mathrm{S}_{2p}$ and thus $G_\Delta\cong G^{\mathcal{C}}_{\Delta}=\mathrm{A}_{2p-1}$ or $\mathrm{S}_{2p-1}$. The kernel of the action of $G_\Delta$ on $\Delta$ is normal in $G_\Delta$ and so can only be $1$, $\mathrm{A}_{2p-1}$ or $\mathrm{S}_{2p-1}$. Since  $G_\Delta^\Delta$ is transitive of degree $2p>2$, this kernel must be trivial. Hence $G_\Delta\cong G_\Delta^\Delta$
is primitive  of degree $2p$ and neither $\mathrm{A}_{2p-1}$ nor $\mathrm{S}_{2p-1}$ has such an action, for $p\in\{5,11\}$.
This contradiction implies that $K\ne 1$.

Since $K\ne 1$ and $K$ is normal in $G$, its orbits are nontrivial blocks of imprimitivity for $G$ in $\P$, and by Claim 1, they must have size  $x=2p$. Hence the set of $K$-orbits in $\P$ is the partition  $\cal C$.
 Since $1\ne \Soc(K)\unlhd G$ it follows that $\Soc(K)^\Delta\ne1$ and hence
$\Soc(K)^\Delta$ contains the socle of $G^{\Delta}_{\Delta}$, which is $2$-transitive on $\Delta$ (see above). Therefore $\Soc(K)^\Delta$ is $2$-transitive, and so also $K^\Delta$ is $2$-transitive. By Burnside's Theorem (see~\cite[Theorem 3.21]{PS}), since $|\Delta|=2p$ is not a prime power,
$G^{\Delta}_{\Delta}$ , $K^\Delta$ and $\Soc(K)^\Delta$  are almost simple with 2-transitive nonabelian simple socle. As mentioned above these 2-transitive groups must have socle $\PSL(2,9)$ or $\mathrm{A}_{10}$ for degree 10, and  $M_{22}$ or $\mathrm{A}_{22}$ for degree 22, and that socle is also 2-transitive on $\Delta$.
\qed

\medskip\noindent
\textit{Claim 3:}\quad The group $K$ is faithful on $\Delta$, so $K$ is almost simple with nonabelian simple socle.

\smallskip\noindent
\textit{Proof of Claim:}\quad Let $\Delta\in\cal C$ and suppose that $A=K_{(\Delta)}\ne1$. Let $F$ denote the set of fixed points of $A$, so $\Delta\subseteq F$. If $\beta\in F$ and $\beta\in\Delta'\in\cal C$, then since $K$ is transitive on $\Delta'$ (Claim 2) and $A\unlhd K$, it follows that $A$ fixes $\Delta'$ pointwise. Thus $A\leq K_{(\Delta')}$, and since $K_{(\Delta)}, K_{(\Delta')}$ are conjugate in $G$ we have $A= K_{(\Delta')}$. Therefore $F$ is a union of parts of $\cal C$.

If $g\in G$, then $A^g$ has fixed point set $F^g$ and $F^g$ is a union of some parts of $\cal C$.
Thus if $F\cap F^g$ contains a point $\beta$ and $\beta\in\Delta'\in\cal C$, then by the previous paragraph $A=   K_{(\Delta')} = A^g$ and so $F=F^g$. It follows that $F$ is a block of imprimitivity for $G$ in $\P$, and $F$ is non-trivial since $A\ne 1$. Thus ${\cal C'}:=\{\ F^g \mid g\in G \}$ is a
non-trivial $G$-invariant partition of $\P$. By Claim 1, $|F|=x$, and since $F$ contains $\Delta$ we conclude that $F=\Delta$. This means that  $A^{\Delta'}\ne 1$ for each $\Delta'\in{\cal C}\setminus\{\Delta\}$, and since $K^{\Delta'}$ is $2$-transitive  (Claim 2), it follows that $A^{\Delta'}$ is transitive. Now choose $\alpha, \beta\in F=\Delta$ and let $B_1, B_2\in\cal B$ be the two blocks containing $\{\alpha,\beta\}$. Then $A\leq G_{\alpha\beta}$, and $G_{\alpha\beta}$ fixes $B_1\cup B_2$ setwise. By Claim 1, there exists $\Delta'\in{\cal C}\setminus\{\Delta\}$ such that $|B_1\cap \Delta'|=\ell=2$, and $|B_2\cap \Delta'|=0$ or $2$. Thus $(B_1\cup B_2)\cap \Delta'$ has size between 2 and 4 and is fixed setwise by $A$. This is a contradiction since $A$
is transitive on $\Delta'$ and $|\Delta'|=2p\geq 10$. Therefore $A=1$ so $K$ is faithful on $\Delta$. By Claim 2, $K\cong K^\Delta$ is almost simple with nonabelian simple socle.  \qed

Since $K$ is 2-transitive of degree $c=2p$, as mentioned above, $K$ has only one conjugacy class of subgroups of index $2p$, and so $K$ has a unique $2$-transitive representation of degree $c$, up to permutational equivalence.  It follows that, for $\alpha\in\Delta$, the stabiliser $K_\alpha$ fixes exactly one point in each part of $\cal C$. Let $\beta$ be another point fixed by $K_\alpha$.
Let
$B_1, B_2\in\cal B$ be the two blocks containing $\{\alpha,\beta\}$. By Claim 1,
$|B_i\cap \Delta|=2$ for each $i$ and hence $K_{\alpha\beta}$ fixes setwise $(B_1\cup B_2)\cap
\Delta$, a set of size 2 or 3. On the other hand $K_{\alpha\beta}=K_\alpha$ since $\beta$ is a fixed point of $K_\alpha$, and by Claim 2, $K$ is 2-transitive on $\Delta$, so the $K_\alpha$-orbits in $\Delta$ have sizes $1, c-1$. This final contradiction completes the proof of Theorem~\ref{Th1}.

\qed

\section{Proof of Theorem \ref{Th2}}\label{sec4}

Our first result in this section proves that the designs arising from
Construction \ref{cons} are all $2$-designs, and inherit certain
symmetry properties from those of the input design.
In particular we show that the designs
coming from projective geometries over
a field of three elements give examples for Theorem \ref{Th2}.

\begin{proposition} \label{exist 1}
Let $\mathcal{S=(P,L)}$ be a $2$-$(v,k,1)$ design, $\ell\in\mathcal{L}$ and $G\leq {\rm Aut}(\mathcal{S})$.
\begin{enumerate}
\item[\rm(i)] Then the design $\D(\mathcal{S})$ given  in Construction~\ref{cons} is a non-symmetric $2$-$(v,k-1,k-2)$ design and $G$ is a subgroup of $\Aut(\D(\mathcal{S}))$;
\item[\rm(ii)] Moreover, if $G$ is flag-transitive on $\mathcal{S}$ and $G_{\ell}$ is $2$-transitive on $\ell$, then $G$ is flag-transitive and point-primitive on $\D(\mathcal{S})$;
\item[\rm(iii)] In particular, if  $\cal S$ is the design of points and lines of the projective space $\PG(n-1,3)$ ($n\geq3$), and $G\geq\PSL(n,3)$, then
$\mathcal{D}(\mathcal{S})$ is a non-symmetric $G$-flag-transitive, $G$-point-primitive $2$-$(v,3,2)$ design.
\end{enumerate}
\end{proposition}

\proof
Let $\mathcal{D}=\mathcal{D}(\mathcal{S})$ with block set
$\mathcal{B}=\{\ell\setminus\{\alpha\}\,\mid\,\ell\in \mathcal{L},\,\alpha\in\ell\}$,
so $\mathcal{D=(P,B)}$.
Let $\alpha,\beta$ be distinct points of $\mathcal{P}$.
Then there exists a unique line $\ell\in \mathcal{L}$,
such that $\alpha,\beta\in \ell$.
As $|\ell|=k$, exactly
$k-2$ blocks of $\mathcal{B}$ contain $\alpha$ and $\beta$.
Thus, $\mathcal{D}$ is a  2-$(v,k-1,k-2)$ design, which is nontrivial provided that $3<k$.
By Lemma \ref{condition 1} applied to $\mathcal{S}$,
$|\mathcal{L}|\geq v$, and since $|\mathcal{B}|=k|\mathcal{L}|>|\mathcal{L}|$
it follows that $\mathcal{D}$ is not symmetric.
Moreover, for all $B=\ell\setminus\{\alpha\}\in \mathcal{B}$ 
and for all $g\in G\leq {\rm Aut}(\mathcal{S})$,
we have $\ell^{g}\in \mathcal{L}$ and $\alpha^{g}\in \ell^{g}$,
and so $B^{g}=(\ell \backslash \{\alpha\})^{g}=\ell^{g}\backslash \{\alpha^{g}\}\in \mathcal{B}$. Thus, $G\leq {\rm Aut}(\D)$ and part (i) is proved.

Now assume that $G$ is flag-transitive on $\mathcal{S}$ and $G_{\ell}$ is
$2$-transitive on $\ell$. Let $\alpha\in \ell$ and $B=\ell\setminus\{\alpha\}$.
From the flag-transitivity of $G$, we know that $G$ acts
primitively on the point set $\mathcal{P}$ by~\cite[Propositions 1--3]{HM},
and $G$ acts transitively on the block set $\mathcal{B}$ of $\D$. Furthermore,
$G_{\ell,\alpha}\leq G_{B}$.
Since $G_{\ell}$ is 2-transitive on $\ell$, $G_{\ell,\alpha}$ is transitive on $B$.
Hence $G_{B}$ is transitive on $B$, and so
$G$ is flag-transitive on $\D$ and part (ii) is proved.

In the special case where  $\cal S$ is the design of points and lines of the projective space $\PG(n-1,3)$ ($n\geq3$), and $H=\PSL(n,3)$, $H$ is flag-transitive on $\cal S$ and $H_\ell$ induces the $2$-transitive group $\PGL(2,3)\cong \mathrm{S}_4$ on $\ell$. Thus part (iii) follows from part (i) and (ii) for any group $G$ such that $H\leq G\leq {\rm Aut}(G)$.
\qed

\subsection{Broad proof strategy and the natural projective action}

In the remainder of the paper we assume the following hypothesis:
\begin{hypothesis}\label{H}
Let $\D=(\P,\cal B)$ be a  $2$-$(v,k,2)$ design admitting a flag-transitive point-primitive group $G$ of automorphisms with socle $X=\PSL(n,q)$ for some $n\geq 3$, where $q=p^f$ with prime $p$ and positive integer $f$.
\end{hypothesis}

Observe that $G\cap {\rm P\Gamma L}(n,q)$ has a natural projective action on a vector space $V$ of dimension $n$ over the field $\mathbb{F}_q$. Consider  a point $\alpha$ of $\D$ and a basis $v_{1},v_{2},\ldots,v_{n}$ of the vector space $V$. Since $G$ is primitive on $\P$, the stabiliser $G_\alpha$ is maximal in $G$, and so by Aschbacher's Theorem~\cite{Asch}(see also~\cite{PB}), $G_\alpha$ lies in
one of the geometric subgroup families $\mathcal{C}_i(1\leq i \leq 8)$,
or in the family $\mathcal{C}_9$ of almost simple subgroups not contained in any of these families.
When investigating the subgroups in the Aschbacher families, we make frequent use of the information
on their structures in~\cite[Chap. 4]{PB}.
We will sometimes use the symbol $\tilde{H}$ to indicate that we are giving the structure of the pre-image of $H$ in the corresponding (semi)linear group.

In the next proposition we treat the case where $\P$ is the point set of the projective space $\PG(n-1,q)$ associated with $V$.

\begin{proposition} \label{4}
Assume Hypothesis~$\ref{H}$, and that $\mathcal{P}$ is the point set of
the projective space $\PG(n-1,q)$, with $G$ acting naturally on $\mathcal{P}$.
Then  either
\begin{enumerate}[(a)]
\item $q=3$, $k=3$, $v=(3^{n}-1)/2$ and $\D=\D(\cal S)$ from Construction~$\ref{cons}$, where $\cal S$ is the design of points and lines of $\PG(n-1,3)$; or
\item $q=2$, $k=4$, $v=2$ and $\D$ is the complement of the Fano plane (that is, blocks are the complements of the lines in $\PG(2,2)$).
\end{enumerate}

\end{proposition}

\proof
Let $\alpha,\beta$ be distinct points. Since $\lambda=2$,
there are exactly two blocks $B_{1}$ and $B_{2}$ containing $\alpha$ and $\beta$.
Moreover, $G_{\alpha\beta}$ fixes $B_{1}\cup B_{2}$ setwise,
so $B_{1}\cup B_{2}$ is a union of $G_{\alpha,\beta}$-orbits.
Let $\ell$ be the unique projective line containing $\alpha$ and $\beta$.
Then $G_{\alpha,\beta}$ is transitive on the $v-(q+1)$ points $\mathcal{P}\backslash\ell$ and on $\ell\backslash\{\alpha,\beta\}.$
Hence, either
\begin{enumerate}
\item[1.] $(B_{1}\cup B_{2})\backslash\{\alpha,\beta\}\supseteq \mathcal{P}\backslash\ell$, or
\item[2.] $B_{1}\cup B_{2}=\ell$.
\end{enumerate}

Suppose first that $(B_{1}\cup B_{2})\backslash\{\alpha,\beta\}\supseteq \mathcal{P}\backslash\ell$.
Then $2k-2\geq |B_{1}\cup B_{2}|\geq 2+v-(q+1)$, that is $k-1\geq (v-q+1)/2$.
Now $r(k-1)=2(v-1)$ (Lemma~\ref{condition 1}) and $v=(q^{n}-1)/(q-1)$,
so that
\begin{equation}\label{Eq11}
r=\frac{2(v-1)}{k-1}\leq\frac{4(v-1)}{v-q+1}=4\cdot \left(1+\frac{q-2}{q^{n-1}+\cdots+q^{2}+2}\right)<8.
\end{equation}
Since $r\geq k$, we have that $k\leq 7$. Now combining this with $k-1\geq (v-q+1)/2$,
we have that $12\geq 2(k-1)\geq q^{n-1}+\cdots+q^{2}+2$.
If $n\geq 4$, then $12\geq q^{n-1}+\cdots+q^{2}+2\geq q^3+q^{2}+2\geq 2^3+2^{2}+2=14$, a contradiction. So $n=3$ and  $12\geq q^{2}+2$, which implies that $q\leq 3$. If $q=3$, then $v=13$, and $6\geq k-1\geq (v-q+1)/2$ implies that $k=7$. Now $r(k-1)=2(v-1)$ implies that $r=4$, contradicting $r\geq k$.
Hence $(n,q)=(3,2)$.
Then $v=7$, $k-1\geq 3$, and $r=2(v-1)/(k-1)\leq4$, and so $r\leq 4\leq k$. Since $r\geq k$, we get that $r=k=4$, and thus
 $b=(vr)/k=7$.
Thus, $\mathcal{D}$ is a symmetric $2$-$(7,4,2)$ design with $X=\PSL(3,2)$.
Since $k=4$, and $G_B$ is transitive on the block $B$,
it follows that $B$ does not contain a line of $\PG(2,2)$.
The only possibility is that $B=\P\backslash \ell'$,
where $\ell'$ is a line of $\PG(2,2)$,
that is, the blocks are complements of the lines of $\PG(2,2)$.
Hence $\D$ is the complement of the Fano projective plane and (b) holds.

Now assume that  $B_{1}\cup B_{2}=\ell$, and every block is contained in a line of the projectice space.
We get $2k-2\geq|B_{1}\cup B_{2}|=q+1$, while $q+1=|\ell|>|B_{i}|=k$.
Hence $q>k-1\geq(q+1)/2>q/2$.

Assume that there are $s$ blocks of $\D$ through $\alpha$
contained in the projective line $\ell$.
Since $G$ acts flag-transitively on the projective space $\PG(n-1,q)$,
for any projective line $\ell'$ and any point $\alpha'\in\ell'$,
there are $s$ blocks containing $\alpha'$
that are contained in $\ell'$.
Since for any two distinct points,
there is a unique projective line containing them,
the sets of blocks on $\alpha$ that are contained in
distinct lines $\ell,~\ell'$ through $\alpha$ are disjoint.
Note that there are $(q^{n-1}-1)/(q-1)$
projective lines through $\alpha$,
so the number of blocks through $\alpha$
is $r=s(q^{n-1}-1)/(q-1)$.

As $r(k-1)=2(v-1)$, it follows that
$s(k-1)(q^{n-1}-1)/(q-1)=2((q^{n}-1)/(q-1)-1)$,
so $s(k-1)=2q$.
Then it follows from $q>k-1> q/2$
that $1>2/s>1/2$, and so $s=3$.
Thus there are 3 blocks through $\alpha$ contained in $\ell$,
and $k-1=2q/3$, so $q=3^{f}$ for some $f$,
and $k=2\cdot3^{f-1}+1$.

Assume that there are $c$ blocks of $\D$ contained in the projective line $\ell$.
Since $G$ acts transitively on the projective lines, for any projective line $\ell'$,
there are $c$ blocks contained in $\ell'$. Now, counting the number of flags $(\gamma,B)$ in two ways, where $\gamma\in \ell$ and $B\subseteq \ell$ for a fixed line $\ell$,
we have that $3(q+1)=ck$, so $3(3^{f}+1)=c(2\cdot 3^{f-1}+1)$, which can be rewritten as $3^{f-1}(9-2c)=c-3$. Suppose
$f\geq2$. Then $3$ divides $c$: when $c=3$, the equation cannot hold, and when $c\geq 6$ the left hand side is negative while the right hand side is positive. Hence $f=1$, $q=3$, $k=3$, and $c=4$.
Therefore, the blocks contained in $\ell$ are all the sets $\ell\backslash\{\gamma\}$,
for $\gamma\in \ell$,
and this implies that $\mathcal{B}=\{\ell\backslash\{\gamma\}\,|\,\ell\in \mathcal{L}, \gamma\in \ell\}$.
Therefore, $\D=\D(\cal S)$ is the design in Construction \ref{cons}, where $\cal S$ is the design of points and lines of $\PG(n-1,3)$.
\qed
\medskip

In what follows, we analyse each of the families $\mathcal{C}_{1}$--$\mathcal{C}_{9}$ for $G_\alpha$.

\subsection{$\mathcal{C}_{1}$-subgroups}\label{sec3.1.1}

In this analysis we repeatedly use the Gaussian binomial coefficient $\gbc{m}{i}_q$ for the number of $i$-spaces in an $m$-dimensional space $\mathbb{F}_q^m$, where $0\leq i\leq m$. A straightforward argument counting bases of $\mathbb{F}_q^m$ and its subspaces shows that, for $i\geq1$,
\[
\gbc{m}{i}_q= \frac{(q^m-1)(q^m-q)\cdots(q^{m}-q^{i-1})}{(q^i-1)(q^i-q)\cdots(q^i-q^{i-1})}
= \frac{\prod_{j=1}^i(q^{m-i+j}-1)}{\prod_{j=1}^i(q^j-1)} = \prod_{j=1}^i \frac{q^{m-i+j}-1}{q^{j}-1}.
\]
We use this equality without further comment. We also use the facts that $\gbc{m}{i}_q=\gbc{m}{m-i}_q$, that the number of complements in
$\mathbb{F}_q^m$ of a given $i$-space is $q^{i(m-i)}$, and hence that the number of decompositions $U\oplus W$ of $\mathbb{F}_q^m$ with $\dim(U)=i$ is $\gbc{m}{i}_q\cdot q^{i(m-i)}$.

\begin{lemma}\label{c1'}
Assume Hypothesis~$\ref{H}$. If the point-stabilizer $G_\alpha\in\mathcal{C}_{1}$,
then $G_\alpha$ is the stabiliser in $G$ of an $i$-space and $G\leq\PGaL(n,q)$.
\end{lemma}
\proof
If $G\leq \PGaL(n,q)$ then $G_\alpha$ is the stabiliser in $G$ of an $i$-space, for some $i$, so assume
that $G\nleq\PGaL(n,q)$. Then $G$ contains a graph automorphism of $\PSL(n,q)$, so in particular $n\geq3$,
and $G_{\alpha}$ stabilizes a pair $\{U,W\}$
of subspaces $U$ and $W$, where $U$ has dimension $i$ and $W$ has dimension $n-i$ with $1\leq i< n/2$.
It follows that $G^*:=G\cap\PGaL(n,q)$ has index $2$ in $G$. Moreover, either $U\subseteq W$ or $U\cap W=0$. 

\textbf{Case~1:} $U \subset W$.

In this case, $v$ is the number $\gbc{n}{n-i}_q$ of $(n-i)$-spaces $W$ in $V$, times the number $\gbc{n-i}{i}_q$ of $i$-spaces $U$ in $W$, so
\begin{align*}
v
&=\gbc{n}{n-i}_q\cdot\gbc{n-i}{i}_q = \prod_{j=1}^{n-i} \frac{q^{n-(n-i)+j}-1}{q^{j}-1}\cdot
\prod_{j=1}^i \frac{q^{(n-i)-i+j}-1}{q^{j}-1}\\
&=\left. \prod_{j=1}^{n-i} (q^{i+j}-1) \right/ {\left(\prod_{j=1}^{n-2i} (q^{j}-1)\cdot \prod_{j=1}^i (q^{j}-1)\right)}\\
&= \prod_{j=1}^{i} \frac{q^{i+j}-1}{q^{j}-1} \cdot \prod_{j=1}^{n-2i} \frac{q^{2i+j}-1}{q^{j}-1}.\\
\end{align*}
Then, using the fact that $q^m-1>q^{m-j}(q^j-1)$, for integers $1\leq j<m$,
\[
v > \prod_{j=1}^i q^{i} \cdot \prod_{j=1}^{n-2i} q^{2i} = q^{i^2 + 2i(n-2i)} = q^{i(2n-3i)}.
\]

Consider the following points of $\D$: $\alpha = \{U,W\}$, where $W=\langle v_{1},v_{2},\ldots, v_{n-i}\rangle$ and $U= \langle v_{1},v_{2},\ldots, v_{i}\rangle$, and $\beta = \{U',W\}$, where $U'= \langle v_{1},v_{2},\ldots, v_{i-1}, v_{i+1}\rangle$. Then the $G^*_\alpha$-orbit $\Delta$ containing $\beta$ consists of all the points $\{U'',W\}$
such that the $i$-space  $U''\subset W$ and $\dim(U\cap U'')=i-1$. Thus the cardinality $|\Delta|$
is the number $\gbc{i}{i-1}_q$ of $(i-1)$-spaces $U\cap U''$ in $U$, times the number $\gbc{n-2i+1}{1}_q-1$ of $1$-spaces in $W/(U\cap U'')$ distinct from $U/(U\cap U'')$.
Therefore, since $\gbc{i}{i-1}_q=\gbc{i}{1}_q$,
\[
|G^*_\alpha:G^*_{\alpha\beta}| = |\Delta| = \gbc{i}{1}_q\cdot \left(\gbc{n-2i+1}{1}_q-1\right) =
\frac{q^i-1}{q-1}\cdot \frac{q(q^{n-2i}-1)}{q-1}.
\]

Note that $G_\alpha$ contains a graph automorphism, and each such graph automorphism interchanges $U$ and $W$, and hence does not leave $\Delta$ invariant. Thus the $G_\alpha$-orbit containing $\beta$ has cardinality $2|\Delta|$ (a subdegree of $G$), so by Lemma~\ref{condition 2}(iv), $r$ divides
\[
4|\Delta|=\frac{4q(q^{n-2i}-1)(q^i-1)}{(q-1)^2}.
\]
Note that $(q^j-1)/(q-1)<2q^{j-1}$ for each integer $j>0$.
It follows that
\[
r\leq \frac{4q(q^{n-2i}-1)(q^i-1)}{(q-1)^2}
<4q\cdot2q^{n-2i-1}\cdot2q^{i-1}=16q^{n-i-1}.
\]
Combining this with $r^2>2v$ and $v>q^{i(2n-3i)}$,
we see that
$16^2q^{2(n-i-1)}>2q^{i(2n-3i)}$,
that is,
\begin{equation}\label{Eq12}
2^7>q^{2(i-1)n-3i^2+2i+2}\geq 2^{2(i-1)n-3i^2+2i+2}.
\end{equation}
Since $n>2i$, it follows that $2(i-1)n-3i^2+2i+2 > 4i(i-1)-3i^2+2i+2 = i^2-2i+2$,
and so $i^2-2i-5<0$, which implies $i\leq 3$.

\textbf{Subcase~1.1:} $i=3$.

Then $n>2i=6$. From \eqref{Eq12} we have $2^7>q^{4n-19}\geq 2^{4n-19}$,
which implies $n\leq 6$, a contradiction.

\textbf{Subcase~1.2:} $i=2$.

Then $n>4$. From \eqref{Eq12} we have $2^7>q^{2n-6}\geq 2^{2n-6}$,
which implies $n=5$ or $6$.
Then $r\mid 4q(q+1)^{n-4}$ (for $n=5$ or $6$) and $v>q^{4n-12}$.
Combining this with $r^2>2v$, we deduce $16q^2(q+1)^{2n-8}>2q^{4n-12}$,
that is, $8(q+1)^{2n-8}>q^{4n-14}$. For $n=6,$ this gives $8(q+1)^4>q^{10}$, which is impossible. Thus $n=5$ and $8(q+1)^2>q^6$, so $q=2$ and $v=5\cdot7\cdot31$. On the one hand $r\mid 24$ and on the other hand the condition $r^2>2v$ implies $r\geq 47$, a contradiction.


\textbf{Subcase~1.3:} $i=1$.

Then $n>2$, $r$ divides $4q(q^{n-2}-1)/(q-1)$, and
\[
v=\frac{(q^{n}-1)(q^{n-1}-1)}{(q-1)^2}.
\]
Combining this with the condition $r\mid 2(v-1)$,
we seee that $r$ divides
\begin{align*}
R:&=\gcd\left(2(v-1),\frac{4q(q^{n-2}-1)}{q-1}\right)\\
  &=2\gcd\left(\frac{(q^{n}-1)(q^{n-1}-1)}{(q-1)^2}-1,\frac{2q(q^{n-2}-1)}{q-1}\right),\\
  &=\frac{2q}{(q-1)^2}\cdot
  \gcd\left(q^{2n-2}-q^{n-1}-q^{n-2}-q+2,2(q-1)(q^{n-2}-1)\right).
\end{align*}
Since
\[
(q^{2n-2}-q^{n-1}-q^{n-2}-q+2)-(q-1)^2=(q^n+q^2-q-1)(q^{n-2}-1)
\]
is divisible by $(q-1)(q^{n-2}-1)$, we see that
\[
\gcd\left(q^{2n-2}-q^{n-1}-q^{n-2}-q+2,(q-1)(q^{n-2}-1)\right)
\text{ divides }
(q-1)^2,
\]
and so
\[
\gcd\left(q^{2n-2}-q^{n-1}-q^{n-2}-q+2,2(q-1)(q^{n-2}-1)\right)
\text{ divides }
2(q-1)^2.
\]
Therefore, $R$ divides
\[
\frac{2q}{(q-1)^2}\cdot2(q-1)^2=4q.
\]
Combining this with $r\mid R$, $r^2>2v$ and $v>q^{2n-3}$,
we deduce $16q^2>2q^{2n-3}$.
Therefore, $8>q^{2n-5}\geq 2^{2n-5}$, which leads to $n=3$, and $q<8$. Note that $v=(q^2+q+1)(q+1)$, so $R=\gcd\left(2(v-1),4q\right)=2\gcd\left(q(q^2+2q+2),2q\right)=2q\gcd\left(q^2+2q+2,2\right)$. When $q$ is odd we see that $R=2q$.
Then $r^2>2v$ leads to
$ 2(q^2+q+1)(q+1)<r^2\leq R^2= 4q^2, $ which is not possible.
Hence $q\in\{2,4\}$ and $R=4q$.

First assume that $q=4$. Then $v=105$ and $R=16$.
Combining this with $r\mid R$ and $r^2>2v$, we conclude that $r=16$.
Then it follows from $r(k-1)=2(v-1)$ and $bk=vr$ that $k=14$ and $b=120$.
Since $G$ is block-transitive, it follows that $X:=\Soc(G)=\PSL(3,4)$ has equal length
orbits on blocks, of length dividing $b=120$. This implies that $X$ has a maximal subgroup of index dividing $120$, and hence by \cite[page 23]{Atlas}, we conclude that
$X$ is primitive on blocks, that the stabiliser $X_B$ of a block $B$ is a maximal $\C_5$-subgroup stabilising an $\mathbb{F}_2$-structure $V_0=\mathbb{F}_2^3 <V$, and $X_B$ has two orbits on 1-spaces, and on $2$-spaces in $V$. An easy computation shows that $X_B$ has precisely four orbits on the point set $\P$, of lengths 14, 14, 21, 56: these are subsets of flags $\{U,W\}$ determined by whether $U\cap V_0$ contains a non-zero vector or not, and whether $W\cap V_0$ is a 2-space of $V_0$ or not. Since $X_B$ preserves the $k=14$ points of $B$, it follows that $B$ is equal to one of the $X_B$-orbits of length 14, so that $X$ acts flag-transitively and point-imprimitively on $\D$, contradicting Theorem~\ref{Th1}. (In fact $G_B$ interchanges the two $X_B$-orbits of length $14$ and so $G_B$ does not leave invariant a point-subset of size 14.)


Thus $q=2$. Then $v=21$ and $R=8$ and $G=\PSL(3,2).2\cong \PGL(2,7)$.
This together with $r\mid R$ and $r^2>2v$ implies $r=8$.
Then we derive from $r(k-1)=2(v-1)$ and $bk=vr$ that $k=6$ and $b=28$.
However, one can see from~\cite[II.1.35]{Handbook} that there is no $2$-$(21,6,2)$ design, a contradiction.
We also checked with \textsc{Magma}  that considering every subgroup of index 28 as a block stabiliser, and each of its orbits of size 6 as a possible block,  the orbit of that block under $G$ does not yield a $2$-design.

\medskip
\textbf{Case~2:} $V=U\oplus W$.

In this case the number $v$ of points is the number $\gbc{n}{i}_q$ of $i$-spaces $U$ of $V$, times the number $q^{i(n-i)}$ of complements $W$ to $U$ in $V$, so
\begin{align*}
v&=q^{i(n-i)} \prod_{j=1}^i \frac{q^{n-i+j}-1}{q^j-1},
\end{align*}
so in particular $p\mid v$, and
by Lemma \ref{bound}(iii),  $r_p$ divides 2.

Note that $q^i-1>q^{i-j}(q^j-1)$, for integers $i>j$.
Thus
\begin{align*}
v&> q^{i(n-i)} \prod_{j=1}^i q^{n-i} =q^{i(n-i)}(q^{n-i})^i=q^{2i(n-i)}.
\end{align*}
We consider the point $\alpha=\{U,W\}$ with $U=\langle v_{1},\dots,v_i\rangle, W=\langle v_{i+1},\ldots, v_{n}\rangle$ and the $G^*_\alpha$-orbit $\Delta$ containing $\beta=\{U', W'\}$ with $U'=\langle v_{1},\ldots, v_{i-1},v_{i+1}\rangle,
W'=\langle v_i,v_{i+2},\ldots, v_{n}\rangle$. Then $\Delta$ consists of all $\{U'',W''\}$ with
$\dim(U''\cap U)=i-1, \dim(W''\cap W)=n-i-1, \dim(U''\cap W)= \dim(W''\cap U)=1$,  so $|\Delta|$ is the number $\gbc{i}{1}_q\cdot q^{i-1}$ of decompositions $U=(U''\cap U)\oplus (W''\cap U)$, times the number
$\gbc{n-i}{1}_q\cdot q^{n-i-1}$ of decompositions $W=(U''\cap W)\oplus (W''\cap W)$.  Thus
\[
|G_{\alpha}^*:G_{\alpha\beta}^*|=|\Delta|= q^{i-1}\frac{q^i-1}{q-1}\cdot q^{n-i-1}\frac{q^{n-i}-1}{q-1} =  q^{n-2}\frac{(q^i-1)(q^{n-i}-1)}{(q-1)^2},
\]
and $G$ has a subdegree $|\Delta|$ or $2|\Delta|$.
By Lemma~\ref{condition 2}(iv), $r$ divides $4|\Delta|$. Since $r_p\mid 2$,
we deduce that $r$ divides $4(q^i-1)(q^{n-i}-1)/(q-1)^2$ (and even $2(q^i-1)(q^{n-i}-1)/(q-1)^2$ if $q$ is even). Let $a=1$ if $q$ is even and $2$ otherwise. Then $r$ divides $2^a(q^i-1)(q^{n-i}-1)/(q-1)^2$

Considering the inequality $r^2>2v>2q^{2i(n-i)}$ and the fact that $(q^j-1)/(q-1)<2q^{j-1}$ for each integer $j>0$, it follows that
\begin{equation}\label{Alice}
q^{2i(n-i)}<2^{2a-1}\frac{(q^i-1)^2(q^{n-i}-1)^2}{(q-1)^4}<2^{2a-1}(2q^{i-1})^2(2q^{n-i-1})^2=2^{2a+3}q^{2n-4}.
\end{equation}
Thus $2^{2a+3}>q^{2n(i-1)-2i^2+4}\geq 2^{2n(i-1)-2i^2+4},$ so (since $n>2i$)
$$
2a-2\geq 2n(i-1)-2i^2>4i(i-1)-2i^2=2i(i-2).
$$
Hence $i=1$ or $2$, and the case $i=2$ only happens if $a=2$, that is if $q$ is odd.

Assume $i=2$, so $q$ is odd. Then  $2\geq 2n(i-1)-2i^2=2n-8$, so $n\leq 5$. On the other hand $n>2i$, so $n=5$.
By \eqref{Alice} $q^{12}<2^7q^{6}$, so $q^{6}<2^7$, a contradiction since $q\geq 3$.
Therefore $i=1$.
In this case $v=q^{n-1}\frac{q^{n}-1}{q-1}$ and we compute that $v-1=\frac{q^{n-1}-1}{q-1}\cdot (q^n+q-1)$. Since $r\mid 2(v-1)$, $r$ divides $\gcd(2\frac{q^{n-1}-1}{q-1}\cdot (q^n+q-1),4\frac{q^{n-1}-1}{q-1})=2\frac{q^{n-1}-1}{q-1}\gcd(q^n+q-1,2)=2\frac{q^{n-1}-1}{q-1}$. In other words $a=1$ in the computation above whether $q$ is odd or even.
Then  by \eqref{Alice} $q^{2(n-1)}<2(q^{n-1}-1)^2/(q-1)^2<2(2q^{n-2})^2=2^3q^{2n-4}$, which can be rewritten as
$q^2<2^{3}$, so $q=2$. Thus $v= 2^{n-1}(2^{n}-1)$ and $r$ divides $2(2^{n-1}-1)$, so $r^2>2v$ implies that $ 2^{2n-1}-2^{n-1}<2(2^{n-1}-1)^2=2^{2n-1}-2^{n+1}+2$, which is impossible.
\qed

\begin{lemma}\label{c1}
Assume Hypothesis~\ref{H}, and that  the point-stabilizer $G_\alpha\in\mathcal{C}_{1}$.
Then either
\begin{enumerate}[(a)]
\item $\D=\D(\cal S)$ is as in Construction~\ref{cons}, where $\cal S$ is the design of points and lines of $\PG(n-1,3)$; or
\item $\D$ is the complement of the Fano plane.
\end{enumerate}
\end{lemma}

\proof
By Lemma~\ref{c1'}, $G\leq\PGaL(n,q)$, and
$G_{\alpha}\cong {\rm P_{i}}$ is the stabiliser of
a subspace $W$ of $V$ of dimension $i$, for some $i$.
As we will work with the action on the underlying space $V$ we will usually consider a linear group  $\tilde{G}$ satisfying $\tilde{X}=\SL(n,q)\leq \tilde{G}\leq\GGL(n,q)$, acting unfaithfully on $\P$ with kernel a subgroup of scalars.
By Proposition~\ref{4} we may assume that $i\geq2$. Also, on applying a graph automorphism that interchanges $i$-spaces and $(n-i)$-spaces (and replacing $\cal D$ by an isomorphic design) we may assume further that $i\leq n/2$.
Then $v$ is the number of $i$-spaces:
 \[
 v=\gbc{n}{i}_q= \prod_{j=1}^i \frac{q^{n-i+j}-1}{q^{j}-1}.
 \]
Using the fact that $q^i-1>q^{i-j}(q^j-1)$, for integers $i>j$, it follows that $v>q^{i(n-i)}$.

Consider the following points of $\D$: $\alpha = W$, where $W= \langle v_{1},v_{2},\ldots, v_{i}\rangle$, and $\beta = W'$, where $W'= \langle v_{1},v_{2},\ldots, v_{i-1}, v_{i+1}\rangle$. Then the $\tilde{G}_\alpha$-orbit $\Delta$ containing $\beta$ consists of all the points $W''$
such that   $\dim(W\cap W'')=i-1$. Thus the cardinality $|\Delta|$
is the number $\gbc{i}{i-1}_q$ of $(i-1)$-spaces $W\cap W''$ in $W$, times the number $\gbc{n-i+1}{1}_q-1$ of $1$-spaces in $V/(W\cap W'')$ distinct from $W/(W\cap W'')$.
Therefore, since $\gbc{i}{i-1}_q=\gbc{i}{1}_q$,
\[
|\Delta|=\gbc{i}{1}_q\cdot \left(\gbc{n-i+1}{1}_q-1\right)=\frac{q(q^{i}-1)(q^{n-i}-1)}{(q-1)^2}.
\]
Since $\tilde{G}$ is flag-transitive, $r$ divides $2|\Delta|$ (by Lemma \ref{condition 2}(iv)).
Combining this with $r^{2}>2v$ (Lemma~\ref{condition 1}(iv))
we have that
\[
\frac{2q^2(q^{i}-1)^2(q^{n-i}-1)^2}{(q-1)^4}>\frac{(q^{n}-1)\cdots(q^{n-i+1}-1)}{(q^{i}-1)\cdots(q-1)}>q^{i(n-i)}.
\]
Since $2q^{j-1}>(q^j-1)/(q-1)$ for all $j\in\mathbb{N}$, it follows that
\begin{equation}\label{Alice2}
q^{i(n-i)}<\frac{2q^2(q^{i}-1)^2(q^{n-i}-1)^2}{(q-1)^4}<2q^2(2q^{i-1})^2(2q^{n-i-1})^2=32q^{2n-2}\leq q^{2n+3}
\end{equation}
Hence
\begin{equation}\label{eq3}
2n+3>i(n-i)
\end{equation}
and so $i^2+3>n(i-2)\geq2i(i-2)$, which implies that $i\leq4$. Note from Lemma~\ref{condition 1} that $r\mid 2(v-1)$. Let $R=2\gcd(|\Delta|,v-1)$. As $r$ divides $2|\Delta|$, it follows that $r$ divides $R$ and hence $r\leq R$.

\textbf{Case~1:} $i=4$.

In this case, we derive from~\eqref{eq3} that $n\leq9$. This together with the restriction $n\geq2i=8$ leads to $n=8$ or $9$. We also deduce from~\eqref{Alice2} that $32q^{2n-2}>q^{4(n-4)}$, that is $32>q^{2n-14}$. First assume that $n=8$. Then $32>q^{2}$, so $q\leq 5$.
We get
\[
|\Delta|=\frac{q(q^4-1)^2}{(q-1)^2}
\]
and
\[
v=\frac{(q^8-1)(q^7-1)(q^6-1)(q^5-1)}{(q^4-1)(q^3-1)(q^2-1)(q-1)}\]


We easily compute that $R=4,6,40,10$ when $q=2,3,4,5$ respectively, in each case contradicting $r^2>2v$, since $r\leq R$.

Next assume $n=9$. Then $32>q^{4}$, so $q=2$. We get
\[
|\Delta|=\frac{q(q^4-1)(q^5-1)}{(q-1)^2}=930
\]
and
\[
v=\frac{(q^9-1)(q^8-1)(q^7-1)(q^6-1)}{(q^4-1)(q^3-1)(q^2-1)(q-1)}=3309747.
\]
Therefore $R=124$, again contradicting $r^2>2v$.

\textbf{Case~2:} $i=3$.

In this case, we derive from~\eqref{eq3} that $n\leq11$. Together with the restriction $n\geq2i=6$ leads to $n\in\{6,7,8,9,10,11\}$.

For $n=6$, $|\Delta|=\frac{q(q^3-1)^2}{(q-1)^2}=q(q^2+q+1)^2$, while \[v-1=\frac{(q^6-1)(q^5-1)(q^4-1)}{(q^3-1)(q^2-1)(q-1)}-1=q(q^8+q^7+2q^6+3q^5+3q^4+3q^3+3q^2+2q+1).\] Thus $R= 2\gcd(|\Delta|,v-1)=2q\gcd((q^2+q+1)^2,q^8+q^7+2q^6+3q^5+3q^4+3q^3+3q^2+2q+1)$. Using the Euclidean algorithm, we easily see that \[\gcd(q^2+q+1,q^8+q^7+2q^6+3q^5+3q^4+3q^3+3q^2+2q+1)=1,\] so $R=2q$, contradicting $r^2>2v$.

For $n=7$, $|\Delta|=\frac{q(q^3-1)(q^4-1)}{(q-1)^2}=q(q^2+q+1)(q+1)(q^2+1)$, while \[v-1=\frac{(q^7-1)(q^6-1)(q^5-1)}{(q^3-1)(q^2-1)(q-1)}-1=q(q^2+1)(q^9+q^8+q^7+2q^6+3q^5+2q^4+2q^3+2q^2+2q+1).\] Thus \[R= 2\gcd(|\Delta|,v-1)=2q(q^2+1)\gcd((q^2+q+1)(q+1),q^9+q^8+q^7+2q^6+3q^5+2q^4+2q^3+2q^2+2q+1).\] Using the Euclidean algorithm, we easily see that \[\gcd(q^2+q+1,q^9+q^8+q^7+2q^6+3q^5+2q^4+2q^3+2q^2+2q+1)=1\] and  \[\gcd(q+1,q^9+q^8+q^7+2q^6+3q^5+2q^4+2q^3+2q^2+2q+1)=1,\] so $R=2q(q^2+1)$, contradicting $r^2>2v$.

Assume now that $8\leq n\leq 11$. We deduce from~\eqref{Alice2} that $32q^{2n-2}>q^{3(n-3)}$, that is $32>q^{n-7}$. So there are only a finite number of cases to consider and we easily check that for all of them, $R^2<2v$, a contradiction.



\textbf{Case~3:} $i=2$.

In this case, the point set is the set of $2$-spaces 
and $n\geq 4$, but the above restrictions on $r$ do not lead easily to contradictions as they do for larger values of
$i$. So we have a different approach.
Recall that  $\tilde{X}=\SL(n,q)\leq \tilde{G}\leq\GGL(n,q)$, acting unfaithfully on $\P$ (with kernel a scalar subgroup of $\tilde{G}$).
First we deal with $n=4$. In this case
\[
v=\frac{(q^4-1)(q^{3}-1)}{(q^2-1)(q-1)} = (q^2+1)(q^2+q+1), \quad |\Delta|=\frac{q(q^{2}-1)^2}{(q-1)^2} = q(q+1)^2
\]
and by Lemmas~\ref{condition 1} and \ref{condition 2}, $r^2>2v$ and $r$ divides
\begin{align*}
2\gcd(v-1, |\Delta|) &=  2\gcd(q^4+q^3+2q^2+q, q(q+1)^2)\\
			& = 2q\gcd(q^3+q^2+2q+1,(q+1)^2)\\
			&= 2q\gcd((q+1)^2(q-1)+3q+2,(q+1)^2)\\
			& = 2q\gcd(3q+2,(q+1)^2) = 2q
\end{align*}
which implies $4q^2\geq r^2 > 2v > q^4$, a contradiction. Thus $n\geq5$.

Let $H:= \tilde{G}\cap\GL(n,q)$. Then setwise stabiliser $H_{\{\alpha,\beta\}}$ of the points $\alpha = W=\langle v_1, v_2\rangle$
and $\beta=W'=\langle v_1,v_3\rangle$, fixes setwise the two blocks $B_1, B_2$ of $\D$ containing $\{\alpha,\beta\}$.
Also $H_{\{\alpha,\beta\}}$ leaves invariant the spaces $Y=W+W'=\langle v_1,v_2,v_3\rangle$ and $Y'=W\cap W'=\langle v_1\rangle$,
induces $\GL(n-3,q)$ on $V/Y$ (since even $\SL(V)\cap (\GL(\la v_1\ra)\times \GL(\la v_4,\dots,v_n\ra))$ induces
$\GL(n-3,q)$ on $V/Y$).
Moreover $H_{\{\alpha,\beta\}}$ is transitive on $V\setminus Y$, and has orbits of lengths $1, 2q, q^2-q$ on the 1-spaces in $Y$.
Since $H_{\{\alpha,\beta\}}\cap \tilde{G}_{B_1}=H_{\{\alpha,\beta\}}\cap H_{B_1}$ is normal of index 1 or 2 in $H_{\{\alpha,\beta\}}$, it follows that
$H_{\{\alpha,\beta\}}\cap H_{B_1}$ also induces at least $\SL(n-3,q)$ on $V/Y$  and is transitive on $V\setminus Y$.
Hence  the only non-zero proper subspaces of $V$ left invariant by  $H_{\{\alpha,\beta\}}\cap H_{B_1}$ are
$Y, W, W', Y'$, and if $q=2, 3$ then possibly also the $q-1$ other 2-spaces of $Y$ containing $Y'$.

We claim that $H_{B_1}$ is irreducible on $V$. Suppose to the contrary that $H_{B_1}$ leaves invariant a nonzero proper subspace
$U$. Then also $H_{\{\alpha,\beta\}}\cap H_{B_1}$ leaves $U$ invariant. We see from the previous paragraph that $U$ must be contained in $Y$. If $U=Y$, then as $\tilde{G}_{B_1}$ is transitive on the set $[B_1]$ of points of $\D$ incident with $B_1$, it follows that all such points must be 2-spaces contained in $Y$. This is impossible since $\dim(Y)=3$, while some block, and hence all blocks, must be incident with a pair of 2-spaces which intersect trivially. Thus $U$ is a proper subspace of $Y$. The only 1-space invariant under $H_{\{\alpha,\beta\}}\cap H_{B_1}$ is $Y'$, and if $U=Y'$ then the same argument would yield that all 2-spaces incident with $B_1$ would contain $Y'$, which is not true since  some block, and hence all blocks, must be incident with a pair of 2-spaces which intersect trivially. Thus $\dim(U)=2$, and $U$ is a 2-space of $Y$ containing $Y'$. Since $H_{B_1}$ does not fix $\alpha$ or $\beta$, it follows that $U\ne W$ or $W'$, and hence $q=2$ or $3$, and $U$ is one of the $q-1$ other 2-spaces
containing $Y'$.  Again, since $\tilde{G}_{B_1}$ is transitive on $[B_1]$, each 2-space $\alpha'\in [B_1]$ intersects $U$ in a 1-space. Let $\gamma=W''$ be a 2-space which intersects $\alpha=W$ trivially, and let $B$ be a block of $\D$ containing $\{\alpha,\gamma\}$. Then $H_B$ leaves invariant a 2-space, say $U'$, and we have shown that both $W\cap U'$ and $W''\cap U'$ have dimension 1, so $U'$ is contained in the 4-space $W\oplus W''$. Now the subgroup induced by $H_{\{\alpha,\gamma\}}$ on $W\oplus W''$ contains $\GL(W)\times\GL(W'')$. The orbit of $U'$ under this group has size $(q+1)^2$. However the group $H_{\{\alpha,\gamma\}}\cap H_{B}$ has index at most $2$ in $H_{\{\alpha,\gamma\}}$ and fixes $U'$, so we have a contradiction.  Thus we conclude that $H_{B_1}$ is irreducible.

The irreducible group $H_{B_1}$ has a subgroup $H_{\{\alpha,\beta\}}\cap H_{B_1}$ inducing at least $\SL(n-3,q)$ on
$V/Y$. We will apply a deep theorem from \cite{NieP} which relies on the presence of various prime divisors of the subgroup order $|H_{B_1}|$.
For $b, e\geq 2$, a primitive
prime divisor (ppd) of $b^e-1$ is a prime $r$ which divides $b^e-1$ but which does not divide $b^i-1$ for any $i<e$. Such ppd's are known to exist
unless either $(b,e)=(2,6)$, or $e=2$ and $b=2^s-1$ for some $s$, (a theorem of Zsigmondy, see \cite[Theorem 2.1]{NieP}). Each ppd $r$ of $b^e-1$ satisfies $r\equiv 1\pmod{e}$, and if $r>e+1$ then $r$ is said to be large; usually $b^e-1$ has a large ppd and the rare exceptions are known explicitly, see
\cite[Theorem 2.2]{NieP}.  Also, if $b=p^f$ for a prime $p$ then each ppd of $p^{fe}-1$ is a ppd of $b^e-1$ (but not conversely) and this type of ppd of $b^e-1$ is called basic. We will apply \cite[Theorem 4.8]{NieP} which, in particular, classifies all subgroups $H_{B_1}$ with the following properties:
\begin{enumerate}
\item for some integer $e$ such that $n/2 < e\leq n-4$, $|H_{B_1}|$ is divisible by a ppd of $q^e-1$ and also by a ppd of $q^{e+1}-1$;
\item for some (not necessarily different) integers $e', e''$ such that $n/2 < e'\leq n-3$ and $n/2 < e''\leq n-3$,  $|H_{B_1}|$ is divisible by a large ppd of $q^{e'}-1$ and a basic ppd of $q^{e''}-1$.
\end{enumerate}
Since $|H_{B_1}|$ is divisible by $|\SL(n-3,q)|$, it is straightforward to check, using \cite[Theorems 2.1 and 2.2]{NieP}, that
$H_{B_1}$ has these
properties whenever either  $n\geq 11$ with arbitrary $q$, or $n\in\{9,10\}$ with $q>2$. In these cases we can apply \cite[Theorem 4.8]{NieP} to the irreducible subgroup $H_{B_1}$ of $\GL(n,q)$. Note that $H$ does not contain $\SL(n,q)$ since it fixes $[B_1]$ setwise; also, since $e, e+1$ differ by 1 and $e+1\leq n-3$, $H$ is not one of the `Extension field examples' from   \cite[Theorem 4.8 (b), see Lemma 4.2]{NieP}, and finally since $n\geq9$ and $e+1\leq n-3$, $H$ is not one of the `Nearly simple examples' from \cite[Theorem 4.8 (c)]{NieP}. Thus we conclude that either $n\in\{9,10\}$ with $q=2$, or $n\in\{5, 6, 7, 8\}$.

Finally we deal with the remaining values of $n$. Since $H_{\{\alpha,\beta\}}\cap H_{B_1}$ has index at most 2 in $H_{\{\alpha,\beta\}}$
it follows that $H_{B_1}$ has a subgroup of the form $[q^{3\times (n-3)}].\SL(n-3,q)$ which is transitive on $V\setminus Y$, and hence $H_{B_1}$
has order divisible by $q^{x}$ with $x=x(n)=3(n-3)+\binom{n-3}{2} = (n-3)(n+2)/2$; also $H_{B_1}$ does not contain $\SL(n,q)$ since it fixes
$[B_1]$ setwise. It follows that $H_{B_1}\cap\SL(n,q)$ is contained in a maximal subgroup of $\SL(n,q)$ which is irreducible (that is, not in class $\mathcal{C}_1$ in \cite{Low}) and has order divisible by $q^{x(n)}$.  A careful check of the possible maximal subgroups in the relevant
tables in \cite{Low}, as listed in Table~\ref{tabc1}, shows that no such subgroup exists.
This completes the proof.
\qed

\begin{table}[h]
\begin{center}
\caption {Tables from \cite{Low} to check for the proof of Lemma~\ref{c1}, Case $i=2$}\label{tabc1}
\vspace{5mm}
\begin{tabular}{c|c|l}
\hline
$n$ &$x(n)$& Tables from \cite{Low} for $n$ \\
\hline
$5$  & $7$  &  Tables 8.18 and 8.19   \\
$6$ & $12$  &  Tables 8.24 and 8.25 \\
$7$ &$18$ & Tables 8.35 and 8.36 \\
$8$ & $25$ &   Tables 8.44 and 8.45  \\
$9$ & $33$ &  Tables 8.54 and 8.55 \\
$10$ & $42$&  Tables 8.60 and 8.61  \\
\hline
\end{tabular}
\end{center}
\end{table}

\subsection{$\mathcal{C}_{2}$-subgroups}\label{sec3.2}

Here $G_{\alpha}$ is a subgroup of type $\GL(m,q)\wr {\rm \mathrm{S}_{t}}$, preserving a decomposition $V=V_{1}\oplus\cdots\oplus V_{t}$
with each $V_{i}$ of the same dimension $m$, where $n=mt$, $t\geq2$. We can think of the pointset of $\D$ as the set of these decompositions (for a fixed $m$ and $t$).
Note that graph automorphisms swap $i$-spaces with $n-i$-spaces, so $G\leq \PGaL(n,q)$ unless $t=2$. When $t=2$ we have to consider that $G$ could contain graph automorphisms, and so could $G_\alpha$.

\begin{lemma}\label{c2}
Assume Hypothesis~\ref{H}. Then the point-stabilizer $G_\alpha\notin\mathcal{C}_{2}$.
\end{lemma}

\proof
 Recall that we denote $\gcd(n,q-1)$ by $d$.
By Lemma \ref{eq2},  $|X|=|\PSL(n,q)|>q^{n^2-2}$,
and by \cite[Proposition 4.2.9]{PB},
\[
v=\frac{|\GL(mt,q)|}{|\GL(m,q)|^tt!}\quad\mbox{so}\quad |X_{\alpha}|=\frac{|X|}{v}=\frac{t!|\GL(m,q)|^t}{d(q-1)}.
\]

\textbf{Case~1:} $m=1$.

Then $n=t\geq 3$, so $\tilde{G}\leq\GGL(n,q)$. Take $\alpha$ as the decomposition $\oplus_{i=1}^n\la e_i\ra$ and $\beta$ as the decomposition
 $\la e_1+e_2\ra\oplus (\oplus_{i=2}^n\la e_i\ra)$. The orbit of $\beta$ under $G_{\alpha}$ consists of the decomposition $\la e_i+\lambda e_j\ra\oplus (\oplus_{\ell\neq i}\la e_\ell\ra)$, which has size $s:=n(n-1)(q-1)$.
Thus by Lemmas~\ref{condition 1}(iv) and~\ref{condition 2}(iv), and Table~\ref{tab1},
 \[
4n^2(n-1)^2(q-1)^2\geq (2s)^2\geq r^2 > 2v =2\frac{|\GL(n,q)|}{(q-1)^n n!}  > 2\frac{q^{n^2}}{4(q-1)^nn!}
\]
so $8n^2(n-1)^2 n! > q^{n^2}/(q-1)^{n+2}>q^{n^2-n-2}$. This implies that either $(n,q)=(5,2)$ or $(4,2)$, or  $n=3$ and $q\leq 5$.

Suppose first that $n=3$. Then $v=q^3(q^2+q+1)(q+1)/6$ and $r\leq 2s=12(q-1)$. Since $r^2>2v$ we conclude that $q=2$ or $3$. In either case $v$ is divisible by $q$, and since $r$ divides $2(v-1)$ (Lemma~\ref{condition 1}), $r$ is not divisible by $4$ if $q=2$, and not divisible by $3$ if $q=3$. Hence $r$ divides $6(q-1)$ if $q=2$, or $4(q-1)$ if $q=3$ (Lemma~\ref{condition 2}), and then $r^2>2v$ leads to a contradiction. Thus $q=2$ and $n$ is 4 or 5.
In either case, $v$ is divisible by $4$, so $4$ does not divide $r$ (Lemma~\ref{bound}). Then, since $r$ divides $2s=2n(n-1)$, we see that $r$ divides $6$ or $10$ for $n=4, 5$ respectively, giving  a contradiction to $r^2> 2v$. Thus we may assume that $m\geq2$.

\textbf{Case~2:} $t=2$.

Next we deal with the case where $G$ may contain a graph automorphism, namely the case $t=2$, so $n=mt\geq 4$, and $G$ acts on decomposition into two subspaces of dimension $m=n/2$.
Let ${\alpha}$ be the decomposition $V_{1}\oplus V_{2}$ where
\[
V_1=\langle v_{1},\ldots, v_{m}\rangle, \quad V_2=\langle v_{m+1},\ldots, v_{2m}\rangle.
\]
Leet $\beta$ be the decomposition $V_1'\oplus V_2'$, where $V_1'= \langle v_{1},\ldots, v_{m-1},v_{m+1}\rangle$ and $V_2'=
\langle v_{m},v_{m+2},\ldots, v_{2m}\rangle$. Let $G^*:=G\cap\PGaL(n,q)$, so $|G:G^*|\leq 2$. Since $G$ is point-primitive, $G$ is point-transitive, and so
$|G_\alpha:G^*_\alpha|=|G:G^*|\leq 2$.

Moreover, let $G^*_{V_1,V_2}$ be the subgroup of $G^*_\alpha$ fixing $V_1$ and $V_2$, so $G^*_{V_1,V_2}$  has index at most $2$ in
$G^*_\alpha$. If $m>2$, then we are in the same situation as in Lemma \ref{c1'} (Case 2) with $i=m=n/2$ and \[
|\beta^{G^*_{V_1,V_2}}|= q^{n-2}\frac{(q^m-1)^2}{(q-1)^2} =q^{2(m-1)}\frac{(q^m-1)^2}{(q-1)^2}.
\]
If $m=2$, then we have double counted (as $G^*_{V_1,V_2}$  does not fix each of the spaces $V_i\cap V_j'$; in fact it contains an element $x:v_1\leftrightarrow v_2, v_3\leftrightarrow v_4$), and $|\beta^{G^*_{V_1,V_2}}|=q^{2(m-1)}\frac{(q^m-1)^2}{2(q-1)^2}$.
In both cases, $|\beta^{G_{\alpha}}|\mid 4q^{2(m-1)}\frac{(q^m-1)^2}{(q-1)^2}.$
By Lemma \ref{condition 2}(iv), $r$ divides $2|\beta^{G_{\alpha}}|$, and hence
\[
r\mid 8q^{2(m-1)}\frac{(q^m-1)^2}{(q-1)^2}.
\]
Note that
\[
v=\frac{|X|}{|X_{\alpha}|}=\frac{q^{m^2}(q^{2m}-1)\cdots(q^{m+1}-1)}{2(q^m-1)\cdots(q-1)}>\frac{q^{2m^2}}{2}
\]
and in particular $p\mid v$. By Lemma \ref{bound}(iii),  $r_p$ divides 2,
and hence $r$ divides
$
\frac{8(q^m-1)^2}{(q-1)^2}.
$
This together with $r^2>2v$ leads to
\begin{equation}\label{Eq20}
\frac{64(q^m-1)^4}{(q-1)^4}>q^{2m^2}.
\end{equation}
It follows that $64\cdot(2q^{m-1})^4>q^{2m^2}$ and so
\[
2^{10}>q^{2(m^2-2m+2)}\geq 2^{2(m^2-2m+2)}.
\]
Hence $10>2(m^2-2m+2)$ and so $m=2$ and $r\mid 8(q+1)^2$. Then we deduce from~\eqref{Eq20} that $64(q+1)^4>q^{8}$, which implies that $q=2$ or $3$.
Assume $q=2$.  Then $r_2\mid 2$, so $r$ divides $2(q+1)^2=18$, contradicting the condition $r^2>2v=560$.
Hence $q=3$, $r\mid 2^7$ and $v=5265$.
Combining this with $r\mid2(v-1)$ we conclude that $r$ divides $2^5$, again contradicting the condition $r^2>2v$. Thus $t\geq3$ and in particular
$n=mt\geq6$ and $G\leq \GGL(n,q)$.

\textbf{Case~3:} $t\geq3$.

Since
$|\GL(m,q)|<q^{m^2}$, we have
\[
|X_{\alpha}|=\frac{t!|\GL(m,q)|^t}{d(q-1)}<\frac{t!q^{n^2/t}}{d(q-1)}.
\]
Combining this with the assertion $|X|<2(df)^2|X_{\alpha}|^3$ from Lemma \ref{bound}(i), we obtain
\[
|X|<\frac{2f^2(t!)^3q^{3n^2/t}}{d(q-1)^3}<2(t!)^3q^{3n^2/t}.
\]
It then follows from $|X|>q^{n^2-2}$ that $q^{n^2-2}<2(t!)^3q^{3n^2/t}$, that is,
\begin{equation}\label{Eq18}
q^{n^2(1-\frac{3}{t})-2}<2(t!)^3.
\end{equation}
Since $n\geq 2t$, we derive from~\eqref{Eq18} that
\begin{equation}\label{Eq19}
2^{4t(t-3)-2}\leq q^{4t(t-3)-2}\leq q^{n^2(1-\frac{3}{t})-2}<2(t!)^3.
\end{equation}
Hence either $t=3$ or $(t,q)=(4,2)$. Consider the latter case. Here~\eqref{Eq18} becomes $2^{n^2/4-2}<2\cdot(4!)^3$ and hence $n\leq8$. As $n\geq2t=8$, we conclude that $n=8$ and $m=2$. However, then $|X|=|\PSL(8,2)|$ and $|X_\alpha|=24|\GL(2,2)|^4$, contradicting the condition $|X|<2(df)^2|X_{\alpha}|^3=2|X_{\alpha}|^3$ from Lemma \ref{bound}(i).

Thus $t=3$, and ${\alpha}$ is a decomposition  $V_{1}\oplus V_{2}\oplus V_3$ with $\dim(V_{1})=\dim(V_{2})=\dim(V_{3})=m=n/3$. Say
\[
V_1=\langle v_{1},\ldots, v_{m}\rangle, \quad V_2=\langle v_{m+1},\ldots, v_{2m}\rangle, \quad V_3=\langle v_{2m+1},\ldots, v_{3m}\rangle.\]
Let $\beta$ be the decomposition $\langle v_{1},\ldots, v_{m-1},v_{m+1}\rangle\oplus\langle v_{m},v_{m+2},\ldots, v_{2m}\rangle\oplus V_3$. Arguing as in Case~2 we find that
$|\beta^{G_{V_1,V_2,V_3}}| =q^{2(m-1)}\frac{(q^m-1)^2}{(q-1)^2}$ if $m\geq3$, or $|\beta^{G_{V_1,V_2,V_3}}| =q^{2(m-1)}\frac{(q^m-1)^2}{2(q-1)^2}$ if $m=2$.
Now $G_{V_1,V_2,V_3}$ has index dividing $6$ in $G_\alpha$, so $|\beta^{G_{\alpha}}|$ divides $6q^{2(m-1)}\frac{(q^m-1)^2}{(q-1)^2}$.
By Lemma \ref{condition 2}(iv), $r$ divides $2|\beta^{G_{\alpha}}|$. Since $v=\frac{|\GL(3m,q)|}{|\GL(m,q)|^3 3!}$, it follows that $p$ divides $v$ and so by Lemma~\ref{bound}, $r_p$ divides $2$, and hence
\[
r\mid 12\frac{(q^m-1)^2}{(q-1)^2},\quad \mbox{so}\quad r^2 < 144 (2q^{m-1})^4 = 2304 q^{4m-4}.
\]
Note that
\[
v=\frac{|\GL(3m,q)|}{|\GL(m,q)|^3 3!}=\frac{q^{3m^2}}{6} \prod_{i=1}^m\frac{q^{2m+i}-1}{q^i-1}\cdot \prod_{i=1}^m \frac{q^{m+i}-1}{q^i-1} >\frac{1}{6}q^{3m^2+2m\cdot m+m\cdot m}=\frac{q^{6m^2}}{6},
\]
and since $r^2>2v$, we get
\[
\frac{q^{6m^2}}{3}<2v<r^2< 2304 q^{4m-4},
\]
and so $6912>q^{6m^2-4m+4}\geq 2^{6m^2-4m+4}\geq 2^{20}$,
a contradiction.
\qed

\subsection{$\mathcal{C}_{3}$-subgroups}\label{sec3.3}

Here $G_{\alpha}$ is an extension field subgroup.

\begin{lemma}\label{c3}
Assume Hypothesis~\ref{H}. Then the point-stabilizer $G_\alpha\notin\mathcal{C}_{3}$.
\end{lemma}

\proof
By Lemma \ref{eq2} we have $|X|>q^{n^2-2}$,
and by \cite[Proposition 4.3.6]{PB},
\[
X_{\alpha}\cong\mathbb{Z}_a.\PSL(n/s,q^s).\mathbb{Z}_b.\mathbb{Z}_s,
\]
where $s$ is a prime divisor of $n$, $d=\gcd(n,q-1)$,  $a=\gcd(n/s,q-1)(q^s-1)/(d(q-1))$, and $b=\gcd(n/s,q^s-1)/\gcd(n/s,q-1)$.
Thus,
 \[
 |X_{\alpha}|=\frac{s|\GL(n/s,q^s)|}{d(q-1)}.
 \]

\textbf{Case~1:} $n=s$.

Here $n$ is a prime, $|X_{\alpha}|=n(q^n-1)/(d(q-1))$,
and by Lemma \ref{bound}(i),
\[
|X|<2(df)^2|X_{\alpha}|^3
=\frac{2f^2n^3}{d}\left(\frac{q^n-1}{q-1}\right)^3
<2q^2n^3\cdot(2q^{n-1})^3
=16n^3q^{3n-1}.
\]
Combining this with  $|X|>q^{n^2-2}$ we obtain
\begin{equation}\label{Eq20b}
q^{n^2-3n-1}<16n^3,
\end{equation}
and so $2^{n^2-3n-1}<16n^3$, which implies $n\leq 5$.

%

\textbf{Subcase~1.1:} $n=5$.

In this case~\eqref{Eq20b} implies that $q^{9}<16\cdot5^3$, which leads to $q=2$.
However, this means that $|X|=|\PSL(5,2)|$ and
$|X_\alpha|=5\cdot31$,
contradicting the condition $|X|<2(df)^2|X_{\alpha}|^3=2|X_{\alpha}|^3$ from Lemma \ref{bound}(i).

\textbf{Subcase~1.2:} $n=3$.

Then $X=\PSL(3,q)$, $|X_\alpha|=3(q^2+q+1)/d$
and so $v=q^3(q^2-1)(q-1)/3$.
It follows from Lemma \ref{bound}(ii) that $r$ divides $2df|X_{\alpha}|=6f(q^2+q+1)$.
Combining this with $r^2>2v$, we obtain that
$54f^2(q^2+q+1)^2>q^3(q^2-1)(q-1)$, that is,
\[
54f^2>\frac{q^6-q^5-q^4+q^3}{(q^2+q+1)^2}.
\]
This inequality holds only when
\[
q\in\{2,3,4,5,7,8,9,16,32\}.
\]
Let $R=\gcd(6f(q^2+q+1), 2(v-1))$.
Then $r$ is a divisor of $R$.  For each $q$ and $f$ as above,
the possible values of $v$ and $R$ are listed in Table \ref{tab3}.
\begin{table}[h]
\begin{center}
\caption {Possible values of $q$, $v$ and $R$}\label{tab3}
\vspace{5mm}
\begin{tabular}{clllll|lllllllllllllllllll}
\hline
$q$ && $v$ & $R$& & &&&&                       $q$ & $v$ & $R$ \\
\hline
$2$  && $8$     & $14$  &&&&&&                 $8$     & $75264$     & $146$  \\
$3$ && $144$    & $26$ &&&&&&                  $9$     & $155520$    & $182$  \\
$4$ && $960$    & $14$  &&&&&&                 $16$    & $5222400$   & $182$  \\
$5$ &&$4000$    & $186$&&&&&&                  $32$    & $346390528$ & $6342$  \\
$7$ &&$32928$   & $38$&&&&&&  &&\\
\hline
\end{tabular}
\end{center}
\end{table}
Hence the condition $r^2>2v$ implies that $q\in \{2,3,5\}$.

Assume $q=2$. Then $v=8$ and $r$ and divides $14$. From $r(k-1)=2(v-1)$ and $r\geq k\geq 3$ we deduce that $r=7$ and $k=3$, which contradicts the condition that $bk=vr$. Similarly, we have $q\neq 5$ (two cases to check: $(r,k)\in\{(186,44),(93,87)\}.$)
Hence $q=3$. By Table~\ref{tab3}, $v=144$ and $r$ divides $26$. Then from $r(k-1)=2(v-1)$, $bk=vr$ and $r\geq k\geq 3$, we deduce that $r=26$,
$k=12$ and $b=312$. Since $|X_\alpha|=39$, Lemma \ref{condition 2}(ii) implies that $G>X$.  Since $Out(X)$ has size $2$, we must have $G=X.2$ (with graph automorphism). By flag-transitivity, a block stabiliser must have index $312$ and have an orbit of size $12$. We checked with \textsc{Magma}, considering every subgroup of index $312$, and only one has an orbit of size 12 (which is unique), and the orbit of that block under $G$ does not yield a $2$-design.

\textbf{Case~2:} $n\geq2s$.

By Lemma \ref{eq2} we have
\[
|X_\alpha|=\frac{s|\GL(n/s,q^s)|}{d(q-1)}\leq \frac{s(1-q^{-s})(1-q^{-2s})q^{n^2/s}}{d(q-1)}
<\frac{sq^{n^2/s}}{d(q-1)}.
\]
Moreover $|X_\alpha|_p=s_p\cdot q^{n(n-s)/2s}$ and $|X|_p=q^{n(n-1)/2}$.
We deduce that $p$ divides $v=|X:X_\alpha|$, so by Lemma~\ref{bound}(iii), $r_p$ divides $2$, and

\[
|X|<2(df)^2|X_{\alpha}|_{p'}^2|X_{\alpha}|=2(df)^2|X_{\alpha}|^3/|X_{\alpha}|_p^2
<\frac{2f^2s^3q^{(3n^2/s)-n(n-s)/s}}{(s_p)^2d(q-1)^3}
\leq \frac{n^3}4 q^{(2n^2/s) +n}.
\]
For the last inequality, we used that $s\leq n/2$ and $f^2\leq (q-1)^3$. 
Combining this with  $|X|>q^{n^2-2}$ we obtain
\begin{equation}\label{Eq21}
4q^{(1-2/s)n^2-n-2}\leq n^3.
\end{equation}

\textbf{Subcase~2.1:} $s\geq 3$.

Then $n\geq 2s\geq 6$ and \eqref{Eq21} implies that
\[
n^3\geq 4q^{(1-2/s)n^2-n-2}\geq 4q^{(n^2/3)-n-2}\geq 2^{(n^2/3)-n}.
\]
We easily see that this inequality only holds for $n\leq 6$. Therefore $n=2s=6$, and so \eqref{Eq21} implies that $q=2$.
It follows that $X=\PSL(6,2)$ and $|X_\alpha|=3|\GL(2,8)|=2^3\cdot3^3\cdot7^2$, so we can compute $v=|X|/|X_\alpha|=2^{12}\cdot3\cdot5\cdot31$ and $v-1=11\cdot 173149$. We know that $r\mid 2(v-1)$.
By Lemma \ref{bound}(iii), we also know that $r\mid 2df|X_\alpha|_{p'}=2\cdot3^3\cdot7^2$, thus $r\mid 2$, contradicting $r^2>2v$.


\textbf{Subcase~2.2:} $s=2$.

Then $n=2m\geq4$ and $n$ is even,
\[
|X_\alpha|=\frac{2|\GL(n/2,q^2)|}{d(q-1)}= \frac{2q^{n(n-2)/4}(q^n-1)(q^{n-2}-1)\cdots (q^{2}-1)}{d(q-1)},
\]
and
\[
v=\frac{q^{n^2/4}(q^{n-1}-1)(q^{n-3}-1)\cdots(q^{3}-1)(q-1)}{2}.
\]
As we observed, $r_p\mid 2$. Also $v$ is even, and so, from $r(k-1)=2(v-1)$ we deduce that $4\nmid r$.

First assume that $n=4$. Then
\[
|X_\alpha|=\frac{2q^{2}(q^4-1)(q+1)}{d}
\quad\text{and}\quad
v=\frac{q^4(q^3-1)(q-1)}{2}.
\]
By Lemma \ref{bound}(iii),  $r$ divides $2df|X_\alpha|_{p'}$ and hence $r\mid 2f(q^4-1)(q+1)$, which can be rewritten as  $r\mid 2f(q^2+1)(q-1)(q+1)^2$  .
Note that
\[
v-1=\frac{(q+1)(q^7-2q^6+2q^5-3q^4+4q^3-4q^2+4q-4)}{2}+1,
\]
so that $\gcd(v-1,q+1)=1$. Hence, since $r\mid 2(v-1)$, it follows that $\gcd(r,q+1)\mid 2$.
Moreover, it follows from $(q-1)\mid v$ that $\gcd(r,q-1)\mid 2$.
Combining this with $4\nmid r$ and $r\mid 2f(q^4-1)(q+1)$,
we obtain $r\mid 2f(q^2+1)$. Therefore,  using Lemma~\ref{condition 1}(iv),
\[
4f^2(q^2+1)^2\geq r^2>2v=q^4(q^3-1)(q-1).
\]
However, there is no $q=p^f$ satisfying $4f^2(q^2+1)^2>q^4(q^3-1)(q-1)$, a contradiction.

Thus $n\geq 6$. Recall that  $\tilde{X}=\SL(n,q)\leq \tilde{G}\leq\GGL(n,q)$, acting unfaithfully on $\P$ (with kernel a scalar subgroup of $\tilde{G}$).  We regard $V$ as an $m$-dimensional vector space over $\mathbb{F}_{q^2}$
with basis $\{e_1,e_2,\ldots,e_m\}$ and $\tilde{G}_\alpha$ the subgroup of $\tilde{G}$ preserving this vector space structure. Take $w\in\mathbb{F}_{q^2}\backslash\mathbb{F}_q$. Then
\[
V=\langle e_1,e_2,\ldots,e_m\rangle_{\mathbb{F}_{q^2}}
       =\langle e_1,we_1,e_2,we_2,\ldots,e_m,we_m\rangle_{\mathbb{F}_{q}}.
\]
Let
\[
W=\langle e_1,e_2\rangle_{\mathbb{F}_{q^2}}
       =\langle e_1,we_1,e_2,we_2\rangle_{\mathbb{F}_{q}}.
\]
Consider $g\in \SL(n,q)$ defined by
\[
\begin{cases}
e^g_1=e_1,~e^g_2=-e_2,~(we_1)^g=we_2,~(we_2)^g=we_1   &\text{for }1\leq i\leq 2;\\
(e_i)^g=e_i~,(we_i)^g=we_i  &\text{for }3\leq i\leq m.
\end{cases}
\]
Then $g$ does not fix $\alpha$. Let $\beta=\alpha^g$ and let $\tilde{G}_{\alpha,(W)}$ be the subgroup of $\tilde{G}_\alpha$ fixing every vector of $W$.
Note that $W^g=\langle e_1,we_1,-e_2,we_2\rangle_{\mathbb{F}_{q}}=W$ and so $\tilde{G}_{\alpha,(W)}\leq \tilde{G}_{\alpha,\beta}$.
Now $\SL(n,q)_{\alpha,(W)}$ contains $I_4\times\SL(n/2-2,q^2)$, and since this subgroup intersects the scalar subgroup trivially it follows that ${X}_{\alpha,(W)}$  contains a subgroup isomorphic to $\SL(n/2-2,q^2)$ (and so do ${G}_{\alpha,(W)}$, ${G}_{\alpha,\beta}$, and ${X}_{\alpha,\beta}$).
By Lemma~\ref{L:subgroupdiv}, $r$ divides $4df|X_\alpha|/|\SL(\frac{n}{2}-2,q^2)|=8fq^{2n-6}(q^n-1)(q^{n-2}-1)(q+1)$.
Combining this with $r_p\mid 2$ and $4\nmid r$,
we obtain
\begin{equation}\label{Eq22b}
r\mid 2f(q^n-1)(q^{n-2}-1)(q+1).
\end{equation}
Then from $r^2>2v$ and
\[
2v=q^{n^2/4}(q^{n-1}-1)(q^{n-3}-1)\cdots(q^{3}-1)(q-1)
\]
we deduce that
\begin{equation}\label{Eq22}
4f^2(q^n-1)^2(q^{n-2}-1)^2(q+1)^2
>q^{n^2/4}(q^{n-1}-1)(q^{n-3}-1)\cdots(q^3-1)(q-1),
\end{equation}
 and so
\[
4q^2(q^n)^2(q^{n-2})^2(2q)^2
>q^{n^2/4}q^{n-2}q^{n-4}\cdots q^4q^2=q^{(n^2-n)/2}.
\]
Therefore,
\[
2^4q^{4n}>q^{(n^2-n)/2}.
\]
This implies that
\[
2^4>q^{n(n-9)/2}\ge 2^{n(n-9)/2},
\]
and hence $n\le 8$ (since $n$ is even).

Assume that $n=8$. By  \eqref{Eq22} we have that
\[
4f^2(q^{8}-1)^2(q^{6}-1)^2(q+1)^2
>q^{16}(q^{7}-1)(q^{5}-1)(q^3-1)(q-1),\]
and this implies that $q\in\{2,3,4\}$.  By \eqref{Eq22b}, $r$ divides
$
u:=2f(q^{8}-1)(q^{6}-1)(q+1),
$ and hence $r$ divides
 $R:=\gcd(2(v-1),u)$.
However, for each $q\in\{2,3,4\}$,
we find $R^2<2v$, contradicting the fact that $r^2>2v$.

Hence $n=6$, and here $r\mid 2f(q^6-1)(q^4-1)(q+1)$  by  \eqref{Eq22b}, which can be rewritten as
 $r\mid 2f(q^2-q+1)(q^2+1)(q^3-1)(q-1)(q+1)^3$  . Recall that $r\mid 2(v-1)$, and  in this case
$2(v-1)=q^9(q^5-1)(q^3-1)(q-1)-2$, which is congruent to $6$ module $q+1$. Thus
$\gcd(2(v-1),q+1)=6$, and so $\gcd(r,q+1)$ divides 6.
On the other hand, $(q^3-1)(q-1)$ divides $v$, so  $\gcd(r,(q^3-1)(q-1))$ divides 2. Recall that $4\nmid r$.
We conclude that $r\mid 54f(q^2-q+1)(q^2+1)$.
Thus
\[
2916f^2(q^2-q+1)^2(q^2+1)^2\geq r^2>2v=q^9(q^5-1)(q^3-1)(q-1),
\]
which implies that $q=2$.
It then follows that $v=55,552$ and $r\mid 810$.
However, as $r\mid 2(v-1)$, we conclude that $r\mid 6$, contradicting $r^2>2v$.

\qed

\subsection{$\mathcal{C}_{4}$-subgroups}\label{sec3.4}

Here $G_{\alpha}$ stabilises a tensor product $V_1 \otimes V_2$, where $V_1$ has dimension $a$, for some divisor $a$ of $n$, and $V_2$ has dimension $n/a$, with $2\leq a<n/a$, that is $2\leq a<\sqrt{n}$. In particular $n\geq 6$.  Recall that  $d=\gcd(n,q-1)$.

\begin{lemma}\label{c4}
Assume Hypothesis~\ref{H}. Then the point-stabilizer $G_\alpha\notin\mathcal{C}_{4}$.
\end{lemma}

\proof
According to \cite[Proposition~4.4.10]{PB}, we have
\[
|X_{\alpha}|=\frac{\gcd(a,n/a,q-1)}{d}
\cdot|\PGL(a,q)|\cdot|\PGL(n/a,q)|.
\]
By Lemma \ref{eq2},
\[
|X_{\alpha}|\leq|\PGL(a,q)|\cdot|\PGL(n/a,q)|
<\frac{(1-q^{-1})q^{a^2}}{q-1}\cdot
\frac{(1-q^{-1})q^{n^2/a^2}}{q-1}
=q^{a^2+(n^2/a^2)-2}.
\]
Let $f(a)=a^2+{\frac{n^2}{a^2}}-2=(a+\frac{n}{a})^2-2-2n$. This is a decreasing function of $a$ on the interval
$(2,\sqrt{n})$, and hence $f(a)\leq f(2)= (n^2/4)+2$.
Hence $|X_{\alpha}|<q^{a^2+(n^2/a^2)-2}\leq q^{(n^2/4)+2}$.
By Lemma \ref{bound}(i),
\[
|X|<2(df)^2|X_{\alpha}|^3
<2d^2f^2q^{(3n^2/4)+6}<2q^{(3n^2/4)+10}.
\]
Combining this with the fact that $|X|>q^{n^2-2}$ (from Lemma \ref{eq2}),
we obtain
\[
q^{(n^2/4)-12}<2.
\]
Therefore, $n^2/4\leq 12$, which implies that $n=6$,
and hence that $a=2$. Thus
\[
|X_{\alpha}|=\frac{q^4(q^3-1)(q^2-1)^2}{d}
\quad\text{and}\quad
v=q^{11}(q^6-1)(q^5-1)(q^2+1).
\]
Consequently, $p\mid v$ and $v$ is even. By Lemma \ref{bound}(iii),  $r_p$ divides 2,
$4\nmid r$,
and $r$ divides $2df|X_\alpha|_{p'}$ and hence $r$ divides $2f(q^3-1)(q^2-1)^2$. Note that  $(q^3-1)(q+1)\mid q^6-1$ and $q-1\mid q^5-1$, so  $(q^3-1)(q^2-1)$ divides $v$. We conclude that $\gcd(r,(q^3-1)(q^2-1))$ divides $2$.
Hence, $r\mid 2f(q^2-1)$, contradicting the condition~$r^2>2v$.
\qed

\subsection{$\mathcal{C}_{5}$-subgroups}\label{sec3.5}
Here $G_{\alpha}$ is a subfield subgroup of $G$ of type $\GL(n,q_{0})$,
where $q=p^f=q_{0}^s$ for some  prime divisor $s$ of $f$.

\begin{lemma}\label{c5}
Assume Hypothesis~\ref{H}. Then the point-stabilizer $G_\alpha\notin\mathcal{C}_{5}$.
\end{lemma}

\proof
According to \cite[Proposition 4.5.3]{PB},
\[
|X_{\alpha}|\cong\frac{q-1}{d\cdot
\lcm(q_0-1,(q-1)/\gcd(n,q-1))}\PGL(n,q_{0})
\]
and, setting $d_0= \gcd(n, (q-1)/(q_0-1))$ (a divisor of $d$), by Lemma \ref{gcd}(i) we have
\begin{align}\label{Eq23a}
|X_{\alpha}|&=\frac{d_0}{d}\cdot|\PGL(n,q_{0})|
            =\frac{d_0}{d}\cdot q_{0}^{n(n-1)/2}(q_{0}^n-1)(q_{0}^{n-1}-1)\cdots(q_{0}^2-1).
\end{align}
In particular, the $p$-part $|X_\alpha|_p=q_0^{n(n-1)/2}$ is strictly less than $|X|_p=q^{n(n-1)/2}$, so $v=|X|/|X_\alpha|$ is divisible by $p$, and hence, by Lemma \ref{bound}(iii),  $r_p$ divides 2, and $2(df)^2|X_\alpha|^2_{p'}|X_\alpha|>|X|$.
Hence
\[
q^{n^2-2} < |X| < 2d^2f^2 q_{0}^{n(n-1)/2}\cdot \frac{d_0^3}{d^3}\cdot((q_{0}^n-1)
(q_{0}^{n-1}-1)\cdots(q_{0}^2-1))^3.
\]
Since $d_0\leq d<q$, $f<q$ and $2\leq q_0$, this implies that
\begin{equation}\label{Eq23}
q^{n^2-2} < 2d^2f^2\cdot q_{0}^{n(n-1)/2}\cdot q_{0}^{3(n+2)(n-1)/2} < q_0\cdot q^4\cdot q_{0}^{2n^2+n-3}.
\end{equation}
As $q=q_{0}^s$,
we have $s(n^2-2) < 4s+2n^2+n-2$, so
\[
2n^2+n-3\geq s(n^2-6).
\]

\textbf{Case~1:} $s\geq 5$.

Then $2n^2+n-3\geq 5(n^2-6)$, and so $n=3$.
However, the first inequality in \eqref{Eq23} then implies
\[
q^{7}<2\cdot 3^2\cdot q^2\cdot q_{0}^{18},
\]
that is, $q_{0}^{5s-18}<18$. This is not possible as $q_{0}^{5s-18}\geq q_0^7\geq2^7$.

\textbf{Case~2:} $s=3$, that is $q=q_0^3$.

Then $2n^2+n-3\geq 3(n^2-6)$, and so $n=3$ or $4$.
Suppose $n=4$.
Then the first inequality in \eqref{Eq23} implies
\[
q^{14}<2\cdot 4^2\cdot q^2\cdot q_{0}^{33},
\]
that is, $32>q_{0}^3$. This leads to $q_{0}=2$ or $3$, and so $q=q_{0}^3=8$ or $27$,
which does not satisfy the first inequality in \eqref{Eq23}, a contradiction.
Therefore, $n=3=s$, and examining $d=\gcd(3,q-1)$ and $d_0=\gcd(3,q_0^2+q_0+1)$, we see that
$d_0=d\in\{1,3\}$.
The inequality  $|X|<2(df)^2|X_\alpha|^2_{p'}|X_\alpha|$ from Lemma \ref{bound}(iii) becomes (using \eqref{Eq23a})
\[
q^3(q^3-1)(q^2-1)/d=|X|<2d^2f^2q_{0}^3
\cdot(q_{0}^3-1)^3(q_{0}^2-1)^3,
\]
or equivalently, since $q=q_0^3$,
\[
q_{0}^6(q_{0}^9-1)(q_{0}^6-1)<2d^3f^2\cdot(q_{0}^3-1)^3(q_{0}^2-1)^3.
\]
Since $(q_{0}^3-1)^3(q_{0}^2-1)^3<(q_{0}^9-1)(q_{0}^6-1)$ and $d\leq n=3$,
it follows that
\begin{equation}\label{Eq24}
q_{0}^6< 2d^3f^2\leq    54f^2.
\end{equation}
As $3\mid f$ and $q_0=p^{f/3}$,
we then conclude that $f=3$ and $q_{0}=2$, but this means that $d=1$,
contradicting the first inequality of \eqref{Eq24}.

\textbf{Case~3:} $s=2$, that is $q=q_0^2$.

In this case, $d_0=\gcd(n,q_0+1)$ in the expression for $|X_\alpha|$ in \eqref{Eq23a}.
%
Let $a\in \mathbb{F}_q\backslash \mathbb{F}_{q_{0}}$ and consider
$$
g=\begin{pmatrix} a &  &\\  & a^{-1}&\\&&I_{n-2} \end{pmatrix}\in \tilde{X}=\SL(n,q).
$$
Now $g$ does not preserve $\alpha$. Let $\beta=\alpha^g\ne \alpha$. Then
\[
\left\{
\begin{pmatrix}
1 &  &   \\
  &1 &    \\
  &  & B
\end{pmatrix}
\,\middle|\,
B\in \SL(n-2,q_{0})
\right\}
\leq \tilde{X}_{\alpha}\cap (\tilde{X}_{\alpha})^g=\tilde{X}_{\alpha\beta}.
\]
Since this subgroup intersects the scalar subgroup trivially, $X_{\alpha\beta}$ contains a subgroup isomorphic to $\SL(n-2,q_{0})$, and hence so does $G_{\alpha\beta}$.
By Lemma~\ref{L:subgroupdiv},
$r$ divides $4df|X_{\alpha}|/|\SL(n-2,q_{0})|$. Thus, using \eqref{Eq23a},
\[
r\mid 4fd_0q_{0}^{2n-3}(q^n_{0}-1)(q^{n-1}_{0}-1).
\]
Recall that $r_{p}\mid 2$. Moreover,
\[
v=\frac{|X|}{|X_\alpha|}=\frac{q_{0}^{n(n-1)/2}(q_{0}^n+1)(q_{0}^{n-1}+1)\cdots(q_{0}^2+1)}
{d_0}
\]
 is even, and so $4\nmid r$.
Therefore,
\begin{equation}\label{Eq25}
r\mid 2fd_0(q^n_{0}-1)(q^{n-1}_{0}-1).
\end{equation}
From $r^2>2v$, that is to say, $r^2/2 > v$,  we see that
\begin{equation}\label{Eq26}
2f^2d_0^2(q^n_{0}-1)^2(q^{n-1}_{0}-1)^2
>\frac{q_{0}^{n(n-1)/2}(q_{0}^n+1)(q_{0}^{n-1}+1)\cdots(q_{0}^2+1)}{d_0},
\end{equation}
and so, using $f<q=q_0^2$,
\[
2d_0^3q_{0}^{4n+2}>q_{0}^{n^2-1},
\]
that is, $2\gcd(n,q_{0}+1)^3=2d_0^3>q_{0}^{n^2-4n-3}$.
If $n\geq6$, then it follows that $2(q_{0}+1)^3>q_{0}^{9}$, a contradiction.
Thus $3\leq n\leq 5$.

Assume that $n=5$, so $2d_0^3> q_0^2$.
It follows that $d_0\neq 1$, and so $d_0=\gcd(5,q_{0}+1)=5$. This together with $250>q_{0}^2$ implies that $q_{0}\in\{4,9\}$. In either case $f=4$, and the inequality \eqref{Eq26} does not hold, a contradiction. Hence $n\leq 4$.

Since $\PSL(n,q_{0})\lhd X_\alpha$ and $r_{p}\mid 2_{p}$,
by Lemma \ref{parabolic},
$r$ is divisible by the index of a parabolic subgroup of $\PSL(n,q_{0})$, that is, the number of $i$-spaces for some $i\leq n/2$.

\textbf{Subcase~3.1:} $n=4$.
There are $(q_0+1)(q_0^2+1)$ 1-spaces and $(q_0^2+1)(q_0^2+q_0+1)$ 2-spaces, so  $q^2_{0}+1$ divides $r$.
Moreover, it follows from $v=q^6_{0}(q^4_{0}+1)(q^3_{0}+1)(q^2_{0}+1)/\gcd(4,q_{0}+1)$
that $q^2_{0}+1$ divides $v$, since $\gcd(4,q_{0}+1)$ is a divisor of $q_0^3+1$.
Therefore, $q^2_{0}+1$ divides $\gcd(r,v)$.
However $r\mid 2(v-1)$ and hence $\gcd(r,v)\mid 2$,
and this implies that $q^2_{0}+1$ divides $2$, a contradiction.

\textbf{Subcase~3.2:} $n=3$.
Here the number $q_0^2+q_0+1$ of 1-spaces must divide $r$.
Since $r\mid 2(v-1)$ and $q^2_{0}+q_{0}+1$ is odd,
it follows that $q^2_{0}+q_{0}+1$ divides $v-1$.
On the other hand $v=q^3_{0}(q^3_{0}+1)(q^2_{0}+1)/d_0$,
and it follows that $\gcd(v-1,q^2_{0}+q_{0}+1)=q^2_{0}+q_{0}+1$ must divide $2q_{0}+d_0$.
This implies that $q_{0}=2$ and $d_0=\gcd(3,q_0+1)=3$.
Therefore, $7\mid r$, $f=2$ and $v=120$.
However, from \eqref{Eq25} and $r\mid 2(v-1)$ we obtain $r=7$ or $14$, contradicting $r^2>2v$.
\qed

\subsection{$\mathcal{C}_{6}$-subgroups}\label{sec3.6}

Here $G_{\alpha}$ is of type $t^{2m}\cdot\Sp_{2m}(t)$,
where $n=t^m$ for some prime $t\ne p$ and positive integer $m$, and moreover $f$ is odd and is minimal such that $t\gcd(2,t)$ divides $q-1=p^f-1$(see \cite[Table 3.5.A]{PB}).

\begin{lemma}\label{c6}
Assume Hypothesis~\ref{H}. Then the point-stabilizer $G_\alpha\notin\mathcal{C}_{6}$.
\end{lemma}

\proof
From \cite[Propositions~4.6.5 and~4.6.6]{PB} we have $|X_\alpha|\leq t^{2m}|\Sp_{2m}(t)|$,
and from Lemma \ref{eq2} we have $|{\rm Sp}_{2m}(t)|<t^{m(2m+1)}$.
Moreover $t<q$, since $t\gcd(2,t)$ divides $q-1$.
Hence $|X_\alpha|<t^{2m+m(2m+1)}<q^{2m^2+3m}$.
By Lemma \ref{bound}(i), recalling that $d=\gcd(n,q-1)$,
\[
|X|<2(df)^2|X_{\alpha}|^3
<2d^2f^2q^{6m^2+9m}<2q^{6m^2+9m+4}.
\]
Combining this with the fact that $|X|>q^{n^2-2}=q^{t^{2m}-2}$ (by Lemma \ref{eq2}),
we obtain
\[
q^{t^{2m}-(6m^2+9m+6)}<2.
\]
Therefore,
\begin{equation}\label{Eq27}
t^{2m}\leq 6m^2+9m+6.
\end{equation}
As $t\geq 2$, we deduce that $2^{2m}\leq 6m^2+9m+6$, and hence $m\leq 3$.

\textbf{Case~1:} $m=1$.\quad
Here $t=n\geq3$, so $t$ is an odd prime, and from \eqref{Eq27} we have $t^2\leq 21$.
Hence $t=n=3$, so that $t\gcd(2,t)=3$ divides $q-1$, and $d=\gcd(n,q-1)=3$. Also
$|X_\alpha|\leq t^{2m}|\Sp_{2m}(t)|=3^2|{\rm Sp}_2(3)|=2^3\cdot 3^3$, and
then it follows from $q^{n^2-2}<|X|<2(df)^2|X_\alpha|^3$ that
\[
q^7<2(3f)^2|X_\alpha|^3\leq2\cdot(3f)^2\cdot(2^3\cdot 3^3)^3=f^2\cdot2^{10}\cdot 3^{11}.
\]
This inequality, together with the fact that $f$ is odd and is minimal such that  $t\gcd(2,t)=3$ divides $p^f-1$, implies that $q\in\{7,13\}$, and hence also that $f=1$. In particular, $q\equiv4$ or 7$\pmod{9}$, so that, by \cite[Proposition 4.6.5]{PB}, we have $X_\alpha\cong3^2.Q_8$.
According to Lemma \ref{bound}(ii), $r$ divides $2df|X_{\alpha}|=432$. Thus  $r$ divides
$R:=\gcd(432, 2(v-1))$. If $q=7$ then $v= 2^2\cdot7^3\cdot19$, and so $R=6$;  and if $q= 13$, then $v=2^2\cdot7\cdot13^3\cdot 61$, and again $R=6$.
Then $R^2<2v$, contradicting $r^2>2v$.

\textbf{Case~2:} $m=2$.\quad
In this case \eqref{Eq27} shows that $t^4\leq 48$ and so $t=2$ and $n=4$.
Thus $|X_\alpha|\leq t^{2m}|\Sp_{2m}(t)|=2^{4}|{\rm Sp}_{4}(2)|<2^{14}$.
From \cite[Proposition 4.6.6]{PB} we see that $q=p\equiv1\pmod 4$. In particular, $f=1$ and $d=4$.
Then the condition $q^{n^2-2}<|X|<2(df)^2|X_\alpha|^3$ implies that
\[
q^{14}<2\cdot4^2\cdot(2^{14})^3=2^{47},
\]
which yields $q=5$. Then by \cite[Proposition 4.6.6]{PB} we have $X_\alpha\cong2^4.\mathrm{A}_6$.
Therefore, $v=|X|/|X_\alpha|=5^5\cdot13\cdot31$.
By Lemma \ref{bound}(ii),  $r$ divides $2df|X_{\alpha}|=2^{10}\cdot 3^2\cdot5$.
This together with $r\mid 2(v-1)$ implies that $r\mid 4$,
contradicting the condition $r^2>2v$.

\textbf{Case~3:} $m=3$.\quad
We conclude similarly (using \cite[Proposition 4.6.6]{PB}) that $t=2$, $n=8$, $q=p\equiv1\pmod 4$ (so $f=1$) and
$|X_\alpha|<2^{27}$.
However, this together with $q^{n^2-2}<|X|<2(df)^2|X_\alpha|^3$ implies that
$q^{62}<2^{82}\gcd(8,q-1)^2<2^{88}$. Thus $q=2$, a contradiction.
\qed

\subsection{$\mathcal{C}_{7}$-subgroups}\label{sec3.7}

Here $G_\alpha$ is a tensor product subgroup of type $\GL(m,q)\wr {\rm S}_t$,
where $t\geq 2$, $m\geq 3$ and $n=m^t$ (see \cite[Table 3.5.A]{PB}).

\begin{lemma}\label{c7}
Assume Hypothesis~\ref{H}. Then the point-stabilizer $G_\alpha\notin\mathcal{C}_{7}$.
\end{lemma}

\proof
From \cite[Proposition 4.7.3]{PB} we deduce that
$|X_{\alpha}|\leq |\PGL(m,q)|^t\cdot t!$.
This together with Lemma \ref{eq2} implies that $|X_{\alpha}|<q^{t(m^2-1)}\cdot t!$.
Then by Lemma \ref{bound}(i),
\[
|X|<2(df)^2|X_{\alpha}|^3
<2d^2f^2q^{3t(m^2-1)}\cdot (t!)^3<q^{3t(m^2-1)+5}\cdot (t!)^3.
\]
Combining this with the fact that $|X|>q^{n^2-2}=q^{m^{2t}-2}$ (by Lemma \ref{eq2}),
we obtain
\begin{equation}\label{Eq28}
(t!)^3 > q^{m^{2t}-3t(m^2-1)-7} \geq 2^{m^{2t}-3t(m^2-1)-7}.
\end{equation}
Let $f(m)=m^{2t}-3t(m^2-1)-7$. It is straightforward to check that $f(m)$ is an increasing function of $m$, for $m\geq 3$, and hence $f(m)\geq f(3)= 3^{2t}-24t-7$. Thus \eqref{Eq28} implies that
\[
2^{3^{2t}-24t-7}< (t!)^3 \leq t^{3t}.
\]
Taking logarithms to base 2 we have $3^{2t}-24t-7 < 3t\log_2(t)$, which has no solutions for $t\geq2$.
\qed

\subsection{$\mathcal{C}_{8}$-subgroups}\label{sec3.8}

Here $G_\alpha$ is a classical group in its natural representation.

\begin{lemma}\label{c8.1}  Assume Hypothesis~\ref{H}.
 If the point-stabilizer $G_\alpha\in\mathcal{C}_{8}$, then $G_\alpha$ cannot be symplectic.
\end{lemma}

\proof
Suppose for a contradiction that $G_\alpha$ is a symplectic group in $\mathcal{C}_{8}$. Then by \cite[Proposition 4.8.3]{PB}, $n$ is even, $n\geq 4$, and
\[
X_{\alpha}\cong {\rm PSp}(n,q)\cdot\left[\frac{\gcd(2,q-1)\gcd(n/2,q-1)}{d}\right],
\]
where ~$d=\gcd(n,q-1)$. For convenience we will also use the notation $d'= \gcd(n/2,q-1)$ in this proof. Therefore,
\[
|X_\alpha|=q^{n^2/4}(q^n-1)(q^{n-2}-1)\cdots(q^2-1)d'/d,
\]
and so
\[
v=\frac{|X|}{|X_\alpha|}=\frac{q^{(n^2-2n)/4}(q^{n-1}-1)(q^{n-3}-1)\cdots(q^3-1)}{d'},
\]
so in particular $p\mid v$. By Lemma \ref{bound}(iii),  $r_p$ divides 2. 
Since ${\rm PSp}(n,q)\normal X_{\alpha}$, except for  $(n,q)=(4,2)$, we can apply  Lemma \ref{parabolic}, and so in these cases $r$ is divisible by the index of a parabolic subgroup of ${\rm PSp}(n,q)$.
We first treat the case $n=4$.

\textbf{Case~1:} $n=4$.

In this case,
\[
X_{\alpha}\cong {\rm PSp}(4,q)\cdot\left[\frac{\gcd(2,q-1)^2}{\gcd(4,q-1)}\right]
\quad\text{and}\quad
v=\frac{q^2(q^3-1)}{\gcd(2,q-1)}.
\]
If $(n,q)= (4,2)$, then a Magma computation shows that the subdegrees of $G$ are $12$ and $15$, so by
Lemma~\ref{condition 2}~(iv),  $r\mid\gcd(24,30)=6$, contradicting $r^2>2v$. Since $X\cong {\rm A}_8$, using \cite[Theorem 1]{Biplane1} for symmetric designs and \cite[Theorem 1.1]{Liang2} for non-symmetric designs  also rules out this case.
Hence $(n,q)\ne (4,2)$. Then, since the indices of the parabolic subgroups
${\rm P}_1$ and ${\rm P}_2$ in ${\rm PSp}(4,q)$
are both equal to $(q+1)(q^2+1)$, it follows that $(q+1)(q^2+1)\mid r$ and,
since $r\mid 2(v-1)$, that $(q+1)(q^2+1)$ divides $2(v-1)$.
Suppose first that $q$ is even. Then
\[
2(v-1) = 2q^2(q^3-1)-2 = 2(q^2+1)(q^3-q-1) +2q,
\]
which is not divisible by $q^2+1$. Thus $q$ is odd, and we have
\[
2(v-1) = q^2(q^3-1)-2 = (q^2+1)(q^3-q-1) +q-1,
\]
and again this is not divisible by $q^2+1$. Thus $n\ne 4$.

\textbf{Case~2:} $n\geq 6$.

Let $\tilde{X}=\SL(n,q)$, the preimage of $X$ in $\GL(n,q)$, and let  $\{e_1,\dots, e_{n/2}, f_1,\dots, f_{n/2}\}$
be a basis for $V$ such that the nondegenerate alternating form
preserved by $\tilde{X}_\alpha$ satisfies
\[
(e_i,e_j)=(f_i,f_j)=0\quad\mbox{and}\quad (e_i,f_j)=\delta_{ij}\quad\mbox{for all $i, j$}.
\]
Let $\SL(4,q)$ denote the subgroup of  $\tilde{X}$
acting naturally on $U:=\la e_1,e_2, f_1, f_2\ra$ and fixing
$W:=\la e_3,\dots, e_{n/2}, f_3,\dots, f_{n/2}\ra$ pointwise, and let
$\Sp(4,q)=\SL(4,q)\cap \tilde{X}_\alpha$, namely the pointwise stabiliser of $W$ in $\tilde{X}_\alpha$.
Let  $g\in \SL(4,q)\setminus\mathbf{N}_{\SL(4,q)}(\Sp(4,q))$
so $g\not\in\tilde{X}_\alpha$, and let  $\beta=\alpha^g\neq\alpha$. Since $g$ fixes $W$ pointwise,
it follows that the alternating forms preserved by $\alpha$ and $\beta$ agree on  $W$ and hence that
$\tilde{X}_{\alpha\beta}=\tilde{X}_{\alpha}\cap (\tilde{X}_{\alpha})^g$ contains the pointwise stabiliser
$\Sp(n-4,q)$ of $U$ in $\tilde{X}_\alpha$.

Since this subgroup $\Sp(n-4,q)$ intersects the scalar subgroup trivially,  $X_{\alpha\beta}$  contains a subgroup isomorphic to
 $\Sp(n-4,q)$, and hence so does $G_{\alpha\beta}$.
By  Lemma~\ref{L:subgroupdiv},
$r$ divides $4df|X_{\alpha}|/|\Sp(n-4,q)|$,
that is,
\[
r\mid 4d'fq^{2n-4}(q^n-1)(q^{n-2}-1).
\]
Recall that $r_p\mid 2$. Also, since $n\geq6$, $v$ is even,
and hence $4 \nmid r$.
Similarly, it follows from $(q-1)\mid v$ that $r_t\mid2$
for each prime divisor $t$ of $q-1$.
Therefore,
\[
r\mid 2f\cdot\frac{q^n-1}{q-1}\cdot\frac{q^{n-2}-1}{q-1}.
\]
As $f<q$ and $(q^j-1)/(q-1)<2q^{j-1}$ for all $j$, it follows that $r<8q^{2n-3}$.
From $r^2>2v$ we derive that
\begin{align*}
64q^{4n-6}&> 2q^{(n^2-2n)/4}(q^{n-1}-1)(q^{n-3}-1)\cdots(q^3-1)/d'\\
         &>2q^{(n^2-2n)/4}(q^{n-2}q^{n-4}\cdots q^2)/q =2q^{(n^2/2)-n-1},
\end{align*}
and so $32>q^{(n^2/2)-5n+5}\geq 2^{(n^2/2)-5n+5}$,
that is, $n^2-10n<0$. This implies that $n\leq 8$.

Suppose that  $n=8$. Here $d'=\gcd(4,q-1)$. In this case the index of each of the parabolic subgroups
$P_i$, for $1\leq i\leq 4$, is divisible by $q^4+1$, and hence $q^4+1$ divides $r$, which in turn divides $2(v-1)$ by Lemma~\ref{condition 2}.  Then
\[
q^4+1\mid 2d'(v-1) = 2 q^{12}(q^{7}-1)(q^{5}-1)(q^3-1) - 2d'.
\]
Since the remainders on dividing $q^{12}, q^{7}-1, q^{5}-1$ by $q^4+1$ are $-1, -q^3-1$ and $-q-1$, respectively, it follows that
 \[
q^4+1\mid -2 (q^3+1)(q+1)(q^3-1)-2d' = -2 (q^6-1)(q+1)-2d'.
\]
The remainder on dividing $q^6-1$ by $q^4+1$ is $-q^2-1$,
and hence
\[q^4+1 \mid 2 (q^2+1)(q+1)-2d'=2(\frac{q^4-1}{q-1}-d').\]
 This implies that
\[
q^4+1 \mid 2(q^4-1)-2d'(q-1) = 2(q^4+1) -4-2d'(q-1)
\]
and hence $q^4+1\leq 2d'(q-1)+4\leq 8q-4$ (since $d'\leq 4$), a contradiction.

Thus $n=6$. Here $d'=\gcd(3,q-1)$.
The indices of the parabolic subgroups
${\rm P}_1$, ${\rm P}_2$ and ${\rm P}_3$ in ${\rm PSp}(6,q)$ are $(q^3+1)(q^2+q+1)$, $(q^3+1)(q^2+q+1)(q^2+1)$ and $(q^3+1)(q^2+1)(q+1)$, and since one of these numbers divides $r$,
we deduce that $(q^3+1)\mid r$, and so $(q^3+1)$ divides $2d'(v-1)= 2\left(q^{6}(q^{5}-1)(q^3-1) - d'\right)$.
Since the remainders on dividing $q^{6}, q^{5}-1, q^{3}-1$ by $q^3+1$ are $1, -q^2-1$ and $-2$, respectively, it follows that
 $q^3+1$ divides $2 \left( 2(q^2+1)-d'\right)$. Hence
 $q^3+1\leq 2(2q^2+2-d')\leq 2(2q^2+1)$, which implies that $q\leq 4$.
 If $q=3$, but then $d'=1$ and
  $q^3+1=28$ does not divide $2(2(q^2+1)-d')=38$.
 Thus $q$ is even and the divisibility condition implies that $q^3+1$ divides $2(q^2+1)-d'\leq 2q^2+1$, which forces $q=2$ and $d'=1$.
 Hence $v=2^6\cdot7\cdot31$, and therefore $v-1$ is coprime to $5$ and $7$. However $r$, and hence also $2(v-1)$
is divisible by  the index of one of the parabolic subgroups
${\rm P}_1$, ${\rm P}_2$ or ${\rm P}_3$ of ${\rm PSp}(6,2)$, and these are $3^2\cdot7$, $3^2\cdot5\cdot7$, $3^3\cdot5$. This is a contradiction.


\qed

\begin{lemma}\label{c8.2}
Assume Hypothesis~\ref{H}.
If the point-stabilizer $G_\alpha\in\mathcal{C}_{8}$,
then $G_\alpha$ cannot be orthogonal.
\end{lemma}

\proof
Suppose for a contradiction that $G_\alpha$ is an orthogonal group in $\mathcal{C}_{8}$. Then by \cite[Proposition 4.8.4]{PB}, $q$ is odd,
$n\geq3$, and
\[
X_{\alpha}\cong {\rm PSO}^\epsilon(n,q).\gcd(n,2),
\]
where $\epsilon\in\{\circ,+,-\}$. Let $\tilde{X}=\SL(n,q)$ and let $\tilde{X}_\alpha$ denote the full preimage of $X_\alpha$ in $\tilde{X}$.

Let $\varphi$ be the non-degenerate symmetric bilinear form on $V$ preserved by $\tilde{X}_\alpha$, and let $e_1, f_1\in V$
be a hyperbolic pair, that is $e_1, f_1$ are isotropic vectors and $\varphi(e_1, f_1)=1$.
Let $U=\langle e_1,f_1\rangle$, and consider the decomposition $V=U\oplus U^\perp$. Let $g\in \tilde{X}$ fixing $U^\perp$ pointwise  and mapping $e_1$ onto itself and $f_1$ onto $e_1+f_1$. Then $g$ maps the isotropic vector $f_1$ onto the non-isotropic vector  $e_1+f_1$, and so
$g\notin \tilde{X}_{\alpha}$. Let $\beta=\alpha^g$, so that $\tilde{X}_\beta$ leaves invariant the form $\varphi^g$.
Then, since $\varphi$ and $\varphi^g$ restrict to the same form on $U^\perp$, we have that
\[
\left\{
\begin{pmatrix}
I_2   &  \\
    & B
\end{pmatrix}
\,\middle|\,
B\in {\rm SO}^\epsilon(n-2,q)
\right\}
\leq \tilde{X}_{\alpha}\cap \tilde{X}^g_{\alpha}=\tilde{X}_{\alpha\beta}.
\]
Since this group intersects the scalar subgroup trivially,  ${X}_{\alpha\beta}$  contains a subgroup isomorphic to ${\rm SO}^\epsilon(n-2,q)$, and hence so does $G_{\alpha\beta}$.
By Lemma~\ref{L:subgroupdiv},
\begin{align}\label{EqOrth}
r\mid 4df|X_{\alpha}|/| {\rm SO}^\epsilon(n-2,q)|.
\end{align}

We now split into cases where $n$ is odd or even.

\textbf{Case~1:} $n=2m+1$ is odd, so $\epsilon=\circ$ and is usually omitted.

In this case, $X_{\alpha}\cong {\rm PSO}(2m+1,q)$. 
Thus
\[
|X_\alpha|=q^{m^2}(q^{2m}-1)(q^{2m-2}-1)\cdots(q^2-1),
\]
and so
\[
v=|X|/|X_\alpha| = q^{m^2+m}(q^{2m+1}-1)(q^{2m-1}-1)\cdots(q^{3}-1)/d,
\]
where $d=\gcd(2m+1,q-1)$, and this implies that $v$ is even and $p\mid v$.  By Lemma \ref{bound}(iii),  $r_p$ divides 2, so $r_p=1$ since $q$ is odd. Moreover, since $r\mid 2(v-1)$, it follows that $4\nmid r$.

\textbf{Subcase~1.1:} $m=1$.

Then
\[
|X_{\alpha}|=q(q^2-1)
\quad\text{and}\quad
v=q^{2}(q^3-1)/d.
\]
As $p\mid v$, it follows from Lemma \ref {bound}(iii) that
$r$ divides $2df|X_{\alpha}|_{p'}$ and hence
$r$ divides $2df(q^2-1)$.
Combining this with $r\mid 2(v-1)$, we deduce that $r$ divides
\begin{align*}
2\gcd\left(d(v-1),df(q^2-1)\right)
=& 2\gcd\left(q^{2}(q^3-1)-d,df(q^2-1)\right).
\end{align*}
Noting that $\gcd\left(q^{2}(q^3-1)-d,q^2-1\right)$ divides
$$
q^{2}(q^3-1)-d-(q^2-1)(q^3+q-1)=q-1-d,
$$
we conclude that $r$ divides $2df(q-1-d)$.
If $d=\gcd(3,q-1)=3$, then $q\geq7$ (since $q$ is odd) and $r\mid 6f(q-4)$. From $r^2>2v=2q^{2}(q^3-1)/3$ we derive that $54f^2(q-4)^2>q^{2}(q^3-1)$,
which yields  a contradiction. Consequently, $d=1$. Then $r\mid 2f(q-2)$, and from $r^2>2v=2q^{2}(q^3-1)$ we derive that $2f^2(q-2)^2>q^{2}(q^3-1)$, which is not possible.

\textbf{Subcase~1.2:} $m\geq 2$.
By \eqref{EqOrth},
$r\mid 4df|X_{\alpha}|/| {\rm SO}^\epsilon(n-2,q)|$,
that is,
\[
r\mid 4dfq^{2m-1}(q^{2m}-1).
\]
Recall that $r_p=1$ and $4\nmid r$.
We conclude that
\[
r\mid 2df(q^{2m}-1).
\]
Therefore, as $r^2>2v$, we have
\[
4d^2f^2(q^{2m}-1)^2>\frac{2q^{m^2+m}(q^{2m+1}-1)(q^{2m-1}-1)\cdots(q^{3}-1)}{d},
\]
and hence
\begin{align*}
2q^3\cdot q^2\cdot q^{4m}&>2d^3f^2(q^{2m}-1)^2\\
&>q^{m^2+m}(q^{2m+1}-1)(q^{2m-1}-1)\cdots(q^{3}-1)\\
&>q^{m^2+m}(q^{2m}q^{2m-2}\cdots q^{2})\\
&=q^{2m^2+2m},
\end{align*}
This implies that $q^{2m^2-2m-5}<2$ and so $2m^2-2m-5\leq 0$. Thus $m=2$ and $d\leq 5$. Therefore
$q^{6}(q^5-1)(q^{3}-1)<2d^3f^2(q^{4}-1)^2<250f^2(q^{4}-1)^2$, which implies $q=2$, a contradiction.

\textbf{Case~2:} $n=2m$ is even, where $m\geq 2$ since $2m=n\geq 3$.

In this case, $X_{\alpha}\cong {\rm PSO}^\epsilon(2m,q)\cdot 2$ with $\epsilon=\pm$ (we identify $\pm$ with $\pm1$ for superscripts).
Hence
\[
|X_\alpha|=q^{m(m-1)}(q^{m}-\epsilon)(q^{2m-2}-1)(q^{2m-4}-1)\cdots(q^2-1),
\]
and so
\[
v=\frac{|X|}{|X_\alpha|}
 =\frac{q^{m^2}(q^{m}+\epsilon)(q^{2m-1}-1)(q^{2m-3}-1)\cdots(q^{3}-1)}
 {d},
\]
where $d=\gcd(2m,q-1)$, and this implies that $v$ is even and $p\mid v$.  By Lemma \ref{bound}(iii),  $r_p$ divides 2, so $r_p=1$ since $q$ is odd. Moreover, since $r\mid 2(v-1)$, it follows that $4\nmid r$.

By \eqref{EqOrth},
$r\mid 4df|X_{\alpha}|/| {\rm SO}^\epsilon(n-2,q)|$,
that is, $r$ divides
\[4dfq^{2m-2}(q^{m}-\epsilon)\frac{q^{2m-2}-1}{q^{m-1}-\epsilon}=
4dfq^{2m-2}(q^{m}-\epsilon)(q^{m-1}+\epsilon).
\]
As $r_p=1$ and $4\nmid r$, it follows that
\begin{align}\label{Eq30a}
r\mid 2df(q^{m}-\epsilon)(q^{m-1}+\epsilon).
\end{align}
Then we deduce from $r^2>2v$ that
\begin{align}\label{Eq30}
&2d^3f^2(q^{m}-\epsilon)^2(q^{m-1}+\epsilon)^2\nonumber\\
>&q^{m^2}(q^{m}+\epsilon)(q^{2m-1}-1)(q^{2m-3}-1)\cdots(q^{3}-1),
\end{align}
and so
\begin{align*}
2q^3\cdot q^2(2q^{2m-1})^2&>2d^3f^2(q^{m}-\epsilon)^2(q^{m-1}+\epsilon)^2\\
&>q^{m^2}(q^{m}+\epsilon)(q^{2m-1}-1)(q^{2m-3}-1)\cdots(q^{3}-1)\\
&>q^{m^2}(2q^{m-1})(q^{2m-2}\cdots q^{2})\\
&=2q^{2m^2-1}.
\end{align*}
Hence $q^{2m^2-4m-4}<4$ and so $2m^2-4m-4<2$, which implies $m=2$ and $d\leq 4$.
Thus $X_\alpha\cong {\rm PSO}^\epsilon(4,q)\cdot 2$.

Suppose $\epsilon=-$, so that $X_\alpha\cong {\rm PSO}^{-}(4,q)\cdot 2$.
Then \eqref{Eq30} gives
\[2d^3f^2(q^{2}+1)^2(q-1)^2>q^{4}(q^{2}-1)(q^{3}-1),
\]
which can be simplified to
\begin{equation}\label{Eq31bis}
2d^3f^2(q^{2}+1)^2>q^{4}(q+1)(q^{2}+q+1).
\end{equation}
Thus $128f^2(q^{2}+1)^2>q^{4}(q+1)(q^{2}+q+1)$. Since $q$ is odd, this implies that $q=3$ so that $d=2$, but then \eqref{Eq31bis} is not satisfied.

Therefore $\epsilon=+$, so that $X_\alpha\cong {\rm PSO}^{+}(4,q)\cdot 2$.
Then \eqref{Eq30} gives
\begin{equation}\label{Eq31}
2d^3f^2(q^{2}-1)^2(q+1)^2
>q^4(q^2+1)(q^3-1)
\end{equation}
and thus
\[
128f^2(q+1)^2>(q^2+1)(q^3-1).
\]
Since $q$ is odd, we conclude that $q=3$ or  $5$.
However,  $q=3$  does not satisfy \eqref{Eq31}, thus $q=5$, $f=1$ and $d=4$. 
Then $v=|X|/|X_\alpha|=503750$. By \eqref{Eq301}, $r\mid 2df(q^2-1)(q+1)=2^7*3^3$. This together with $r\mid 2(v-1)$ (Lemma \ref{condition 1}(i)) leads to $r\mid 2$, contradicting $r^2>2v$.
\qed

\begin{lemma}\label{c8.3}
Assume Hypothesis~\ref{H}.
If the point-stabilizer $G_\alpha\in\mathcal{C}_{8}$,
then $G_\alpha$ cannot be unitary.
\end{lemma}

\proof
Suppose that $G_\alpha$ is a unitary group in $\mathcal{C}_{8}$. Then by \cite[Proposition 4.8.5]{PB}, $n\geq 3$, $q=q^2_{0}$, and
\[
X_{\alpha}\cong {\rm PSU}(n,q_0)\cdot\left[\frac{\gcd(n,q_0+1)c}{d}\right],
\]
where $d=\gcd(n,q-1)$ and $c=(q-1)/\lcm(q_0+1,(q-1)/d)$.
By Lemma \ref{gcd}(iii), $c=\gcd(n,q_0-1)$.
Hence
\begin{align*}
|X_{\alpha}|
&=|{\rm PSU}(n,q_0)|\cdot\frac{\gcd(n,q_0+1)\gcd(n,q_0-1)}{\gcd(n,q^2_0-1)}\\
&=\frac{c}{d}\cdot q^{n(n-1)/2}_0\prod_{i=2}^n(q^i_0-(-1)^i)
\end{align*}
and
\[
v=\frac{|X|}{|X_\alpha|}
 =\frac{1}{c}\cdot q^{n(n-1)/2}_0\prod_{i=2}^n(q^i_0+(-1)^i),
\]
which implies that $p\mid v$ and $v$ is even.
Since $r\mid 2(v-1)$, it follows that $r_p\mid 2$ and $4\nmid r$.

\textbf{Case~1:} $n=3$.

In this case,
\[
|X_{\alpha}|=\frac{cq^3_0(q^3_0+1)(q^2_0-1)}{d}
\quad\text{and}\quad
v=\frac{q^3_0(q^3_0-1)(q^2_0+1)}{c},
\]
where $c=\gcd(3, q_0-1)$ and $d=\gcd(3,q^2_0-1)$.
Since ${\rm PSU}(n,q_0)\normal X_\alpha$, by Lemma \ref{parabolic}, $r$ is divisible by the index of a parabolic subgroup of ${\rm PSU}(3,q_0)$, that is, $q^3_0+1$.
Hence $(q^3_0+1)\mid r$, which implies that $(q^3_0+1)$ divides $2(v-1)$ and hence also
$2c(v-1)=2 q^3_0(q^3_0-1)(q^2_0+1) -2c$.
Since the remainders on dividing $q^3_0, q^3_0-1$ by $q^3_0+1$ are $-1, -2$, respectively, it follows that
 \[
q_0^3+1\mid 4 (q^2_0+1) - 2c,
\]
which implies that $q_0=2$, $d=3$, $c=1$, and $f=2$.
Thus $v=q^3_0(q^3_0-1)(q^2_0+1)=280$ and $|X_\alpha|=q^3_0(q^3_0+1)(q^2_0-1)/3=72$.
Since $r\mid 2(v-1)$ and $r\mid 2df|X_\alpha|_{p'}$ by Lemma \ref{bound}(iii),
we conclude that $r$ divides $18$, contradicting the condition $r^2>2v$.

\medskip
\textbf{Case~2:} $n\geq 4$.

 Let $\tilde{X}=\SL(n,q)$ and let $\tilde{X}_\alpha$ denote the full preimage of $X_\alpha$ in $\tilde{X}$.
Let $U=\langle e_1, f_1\rangle$ be a nondegenerate $2$-subspace of $V$ relative to the
unitary form $\varphi$ preserved by $\tilde{X}_\alpha$.
Let $A\in \SL(U)$ such that $A$ does not preserve modulo scalars the restriction of $\varphi$ to $U$.
Then the element
$g=\begin{pmatrix}
A &  \\
  & I
\end{pmatrix}\in \tilde{X}$ but $g$ does not lie in $\tilde{X}_\alpha$. Hence  $\beta :=\alpha^g\neq\alpha$.
On the other hand
\[
\left\{
\begin{pmatrix}
I &  \\
  &  B
\end{pmatrix}
\,\middle|\,
B\in {\rm SU}(n-2,q_0)
\right\}\leq \tilde{X}_{\alpha}\cap \tilde{X}^g_{\alpha}=\tilde{X}_{\alpha\beta}.
\]
Since this group intersects the scalar subgroup trivially,  ${X}_{\alpha\beta}$  contains a subgroup isomorphic to ${\rm SU}(n-2,q)$, and hence so does $G_{\alpha\beta}$.
By Lemma~\ref{L:subgroupdiv},
$r$ divides $4df|X_{\alpha}|/|{\rm SU}(n-2,q_0)|$,
that is,
\[
 r\mid 4cfq^{2n-3}_0(q^{n}_0-(-1)^n)(q^{n-1}_0-(-1)^{n-1}).
\]
Since $r_p\mid 2$ and $4\nmid r$, we derive that
\[
r\mid 2cf(q^{n}_0-(-1)^n)(q^{n-1}_0-(-1)^{n-1}).
\]
This together with $r^2>2v$ and $v=|X|/|X_\alpha|$ leads to $r^2|X_\alpha|>2|X|$.
By Lemma \ref{eq2} we have
\[
|X|>q^{2n^2-4}_0
\quad\text{and}\quad
|X_\alpha|<\frac{q^{n^2-1}_0c\gcd(n,q_0+1)}{d}.
\]
Consequently, noting that $\gcd(n,q_0+1)\leq d=\gcd(n,q^2_0-1)$, $c=\gcd(n,q_0-1)<q_0$, and $f<q=q_0^2$, we get
\begin{align*}
2q^{2n^2-4}_0
&<4c^3f^2(q^{n}_0-(-1)^n)^2(q^{n-1}_0-(-1)^{n-1})^2\cdot \frac{q^{n^2-1}_0\gcd(n,q_0+1)}{d}\\
&<4q^7_0(q^{n}_0-(-1)^n)^2(q^{n-1}_0-(-1)^{n-1})^2\cdot q^{n^2-1}_0\\
&<4q^{n^2+6}_0(2q^{n+n-1}_0)^2=16q^{n^2+4n+4}_0
\end{align*}
and hence
\[
q^{n^2-4n-8}_0<8.
\]
It follows that $n^2-4n-8<3$, which implies $n=4$ or $5$.

\textbf{Subcase~2.1:} $n=4$.

Then
\[
v=q^6_0(q^4_0+1)(q^3_0-1)(q^2_0+1)/c,
\]
where $c=\gcd(4,q_0-1)$.
Since $r$ is divisible by the index of a parabolic subgroup of ${\rm PSU}(4,q_0)$,
which is either $(q^2_0+1)(q^3_0+1)$ or $(q_0+1)(q^3_0+1)$,
we derive that $(q^3_0+1)\mid r$. Hence $(q^3_0+1)$ divides $2(v-1)$, and hence also
$2c(v-1)= 2 q^6_0(q^4_0+1)(q^3_0-1)(q^2_0+1) -2c$.
Since the remainders on dividing $q^{6}_0, q^{4}_0+1, q^{3}_0-1$ by $q^3_0+1$ are $1, -q_0+1$ and $-2$, respectively, it follows that
 $q^3_0+1$ divides $2 (-q_0+1)(-2)(q_0^2+1) -2c = 4(q_0-1)(q_0^2+1)-2c$, which equals
 $4(q_0^3+1) - 4(q_0^2-q_0+2) -2c$. It follows that
 $q^3_0+1$ divides $4(q^2_0-q_0+2)+2c$, which implies $q_0=2$.
 Thus $v=2^6\cdot5\cdot7\cdot17$, and the index of a parabolic subgroup of ${\rm PSU}(4,q_0)$ is either $45$ or $27$.
However, neither $45$ nor $27$ divides $2(v-1)$, a contradiction.

\textbf{Subcase~2.2:} $n=5$.

Then
\[
v=q^{10}_0(q^5_0-1)(q^4_0+1)(q^3_0-1)(q^2_0+1)/c,
\]
where $c=\gcd(5,q_0-1)$.
Since $r$ is divisible by the index of a parabolic subgroup of ${\rm PSU}(5,q_0)$,
which is either $(q^2_0+1)(q^5_0+1)$ or $(q^3_0+1)(q^5_0+1)$,
we derive that $(q^5_0+1)\mid r$. Hence $(q^5_0+1)$ divides $2(v-1)$, and hence also
$2c(v-1)= 2q^{10}_0(q^5_0-1)(q^4_0+1)(q^3_0-1)(q^2_0+1) -2c$.
Since the remainders on dividing $q^{10}_0, q^{5}_0-1, (q^{3}_0-1)(q_0^2+1)$ by $q^5_0+1$ are $1, -2$ and $q_0^3-q_0^2-2$,
respectively, it follows that
 $q^5_0+1$ divides $-4(q_0^4+1)(q_0^3-q_0^2-2)-2c$, which equals
 \[
 -4(q_0^5+1)(q_0^2-q_0) -4(-2q_0^4+q_0^3-2q_0^2+q_0-2) -2c.
 \]
 Thus $q^5_0+1$ divides $8q_0^4 -4q_0^3 +8q^2_0 -4q_0+8-2c$.
 However, there is no prime power $q_0$ satisfying this condition, a contradiction.
\qed

\subsection{$\mathcal{C}_{9}$-subgroups}\label{sec3.9}

Here $G_\alpha$ is an almost simple group not contained in any of the subgroups in $\mathcal{C}_1$--$\mathcal{C}_8$.
\begin{lemma}\label{c9}
Assume Hypothesis~\ref{H}.
Then the point-stabilizer $G_\alpha\notin\mathcal{C}_{9}$.
\end{lemma}

\proof
By Lemma \ref{condition 2}(i) and Lemma \ref{eq2},
we have $|G_\alpha|^3>|G|\geq|X|=|\PSL(n,q)|>q^{n^2-2}$.
Moreover, by \cite[Theorem 4.1]{Liemaximal}, we have that $|G_\alpha|<q^{3n}$.
Hence $q^{n^2-2}<|G_\alpha|^3<q^{9n}$, which yields $n^2-2<9n$ and so $3\leq n\leq9$.
Further, it follows from \cite[Corollary 4.3]{Liemaximal} that either $n=y(y-1)/2$ for some integer $y$ or $|G_\alpha|<q^{2n+4}$.
If $n=y(y-1)/2$, then as $3\leq n\leq9$ we have $n=3$ or $6$.
If $|G_\alpha|<q^{2n+4}$, then we deduce from $|G_\alpha|^3>q^{n^2-2}$ that $q^{6n+12}>q^{n^2-2}$,
which implies $6n+12>n^2-2$ and so $3\leq n\leq7$.
Therefore, we always have $3\leq n\leq7$. The possibilities for $X_\alpha$ can be read off from \cite[Tables 8.4, 8.9, 8.19, 8.25, 8.36]{Low}. In Table \ref{tab5} we list all possibilities, sometimes fusing some cases together. Not all conditions from \cite{Low} are listed, but we list what is necessary for our proof. 
Note that in some listed cases  $X_\alpha$ is not maximal in $X$ but there is a group $G$ with $X<G\leq \Aut(X)$  such that $G_\alpha$ is maximal in $G$ and $G_\alpha\cap X$ is equal to this non-maximal subgroup $X_\alpha$.

\begin{longtable}{cclll}
\caption{Possible groups $X$ and $X_\alpha$ }\label{tab5}\\ \hline
\endfirsthead
\multicolumn{5}{l}{}\\
\hline
\endhead
\hline
\multicolumn{5}{r}{}\\
\endfoot\hline
\endlastfoot
Case &$X$         &$X_\alpha$ &      Conditions on $q$ from \cite{Low} &Bound \eqref{lastineq}                          \\
\hline
1&$\PSL(3,q)$  &$\PSL(2,7)$                 &$q=p\equiv1,2,4\pmod7$, $q\neq2$   &$q<14$     \\
2&             &${\rm A}_6$                 &$q=p\equiv1,4\pmod{15}$     &$q<19$             \\
3&             &${\rm A}_6$&$q=p^2,p\equiv2,3\pmod5$, $p\neq3$ &$q<23$\\
\hline
4&$\PSL(4,q)$  &$\PSL(2,7)$            &$q=p\equiv1,2,4\pmod7$, $q\neq 2$  &$q< 4$       \\
5&             &${\rm A}_7$            &$q=p\equiv1,2,4\pmod 7$ &$q< 7$        \\
6&             &${\rm PSU}(4,2)$       &$q=p\equiv1\pmod6$      &$q< 12$                 \\

\hline
7&$\PSL(5,q)$  &$\PSL(2,11)$                &$q=p$ odd &$q< 3$  \\
8&            &${\rm M}_{11}$              &$q=3$    &$q< 4$                               \\
9&            &${\rm PSU}(4,2)$            &$q=p\equiv1\pmod6$  &$q< 5$                     \\
\hline
10&$\PSL(6,q)$ &${\rm A}_6.2_3$             &$q=p$ odd  &$q< 3$                    \\
11&            &${\rm A}_6$                 &$q=p$ or $p^2$ odd                &$q< 2$      \\
12&            &$\PSL(2,11)$           &$q=p$ odd &$q< 3$  \\
13&            &${\rm A}_7$                 &$q=p$ or $p^2$ odd                &$q< 3$              \\
14&            &$\PSL(3,4)^{.}2^-_1$        &$q=p$ odd    &$q< 3$            \\
15&            &$\PSL(3,4)$                 &$q=p$ odd    &$q< 3$             \\
16&            &${\rm M}_{12}$              &$q=3$                        &$q< 4$         \\
17&            &${\rm PSU}(4,3)^{.}2^-_2$   &$q=p\equiv1\pmod{12}$      &$q< 5$              \\
18&            &${\rm PSU}(4,3)$            &$q=p\equiv7\pmod{12}$      &$q< 5$              \\
19&            &$\PSL(3,q)$                  &$q$ odd               &     \\
\hline
20&$\PSL(7,q)$ &${\rm PSU}(3,3)$            &$q=p$ odd        &$q< 2$            \\
\end{longtable}

By Lemma \ref{bound}(i) and  Lemma \ref{eq2}, we have $2d^2f^2|X_\alpha|^3>|X|>q^{n^2-2}$. Using the fact that $d=\gcd(n,q-1)\leq n$, it follows that
\begin{equation}\label{lastineq}
q<\left(2n^2f^2|X_\alpha|^3\right)^{1/(n^2-2)}.
\end{equation}
Note that, except for case (19), we know that $f=1$ or $2$. This inequality gives us, in each case
except (19), an upper bound for $q$, which is listed in the last column in Table \ref{tab5}. Comparing the last two columns of the table we see the condition and bound are satisfied only in the following cases: (1) for $q=11$,
(3) for $q=4$, (5) for $q=2$, (6) for $q=7$, (8) and (16). For case (19), we know that $f<q$ and $|\PSL(3,q)|<q^8$ by Lemma \ref{eq2}, so  $72q^2q^{24}>q^{34}$, that is $q^8<72$, which is not satisfied for any $q$.
In case (1) for $q=11$, $d=1$ and the inequality $2d^2f^2|X_\alpha|^3>q^{n^2-2}$ is not satisfied.

For each of the remaining cases, we compute $v$ and $2df|X_\alpha|$. By Lemma  \ref{bound}(ii), $r\mid 2df|X_\alpha|$. On the other hand $r\mid 2(v-1)$, so  $r$ divides $R:=\gcd(2(v-1),2df|X_\alpha|)$. Now using $R^2\geq r^2>2v$, this argument rules out cases (3) for $q=4$, (6) for $q=7$, (8) and (16). This leaves the single remaining case (5) with $q=2$. Then this argument yields $r\mid 14$, $v=8$.
As $r^2>2v$, $r=7$ or $14$.
 By Lemma \ref{condition 1}(i), $r(k-1)=14$, so the condition $k\geq3$ implies that $r=7$ and $k=3$. Now  Lemma \ref{condition 1}(ii) yields a contradiction since $k\nmid vr$.
Hence, we rule out case (5) for $q=2$, completing the proof.
\qed

\Addresses

\end{document}